\documentclass{article}

\usepackage{arxiv}

\usepackage[utf8]{inputenc} % allow utf-8 input
\usepackage[T1]{fontenc}    % use 8-bit T1 fonts
\usepackage{hyperref}       % hyperlinks
\usepackage{url}            % simple URL typesetting
\usepackage{booktabs}       % professional-quality tables
\usepackage{amsfonts}       % blackboard math symbols
\usepackage{nicefrac}       % compact symbols for 1/2, etc.
\usepackage{microtype}      % microtypography
\usepackage{lipsum}
\usepackage{graphicx}
\graphicspath{ {./images/} }
\usepackage{lipsum}
\usepackage{amsfonts}
\usepackage{graphicx}
\usepackage{epstopdf}
\usepackage{mathtools}
\usepackage{algorithmic}
\usepackage{amsmath}

\newtheorem{assumption}{Assumption}
\newtheorem{theorem}{Theorem}
\newtheorem{lemma}[theorem]{Lemma}

\newtheorem{example}{Example}

\usepackage{chngcntr}
\counterwithin{equation}{section}
\counterwithin{example}{section}
\counterwithin{theorem}{section}
\counterwithin{figure}{section}
\counterwithin{equation}{section}

\usepackage{subfig}
\ifpdf
  \DeclareGraphicsExtensions{.eps,.pdf,.png,.jpg}
\else
  \DeclareGraphicsExtensions{.eps}
\fi

\title{Nonlinear Gradient Mappings and Stochastic Optimization: A General Framework with Applications to Heavy-tail Noise}

\author{
 Du\u{s}an Jakoveti\'{c}\\
  Faculty of Sciences \\
  University of Novi Sad\\
  \texttt{dusan.jakovetic@dmi.uns.ac.rs} \\
   \And
 Dragana Bajovi\'{c} \\
  Faculty of Technical Sciences\\
  University of Novi Sad\\
  \texttt{dbajovic@uns.ac.rs} \\
    \And
 Anit Kumar Sahu \\
  Amazon Alexa AI\\
  \texttt{anit.sahu@gmail.com} \\
    \And
 Soummya Kar \\
  Carnegie Mellon University\\
  \texttt{soummyak@andrew.cmu.edu} \\
    \And
 Nemanja Milo\u{s}evi\'{c} \\
  Faculty of Sciences\\
  University of Novi Sad\\
  \texttt{nmilosev@dmi.uns.ac.rs} \\
    \And
 Du\u{s}an Stamenkovi\'{c} \\
  Faculty of Sciences\\
  University of Novi Sad\\
  \texttt{dusan.stamenkovic@dmi.uns.ac.rs} \\
}

\begin{document}
\maketitle
\begin{abstract}
 We introduce a general framework for nonlinear stochastic gradient descent (SGD) for the scenarios when  gradient noise exhibits heavy tails. 
 The proposed framework subsumes several popular nonlinearity choices, like clipped, normalized, signed or quantized gradient, but we also consider novel nonlinearity choices. We establish for the considered class of methods strong convergence guarantees assuming a strongly convex cost function with Lipschitz continuous gradients under very general assumptions on the gradient noise. Most notably, we show that, for a nonlinearity with bounded outputs and for the gradient noise that may not have finite moments of order greater than one, the nonlinear SGD's mean squared error (MSE), or equivalently, the expected cost function's optimality gap,  converges to zero at rate~$O(1/t^\zeta)$, $\zeta \in (0,1)$. In contrast, for the same noise setting, the linear SGD generates a sequence with unbounded variances. Furthermore, for the nonlinearities that can be decoupled component wise, like, e.g., sign gradient or component-wise clipping, we show that the nonlinear SGD asymptotically (locally) achieves a $O(1/t)$ rate in the weak convergence sense and explicitly quantify the corresponding asymptotic variance. Experiments show that, while our framework is more general than existing studies of SGD under heavy-tail noise, several easy-to-implement nonlinearities from our framework are competitive with state of the art alternatives on real data sets with heavy tail noises.
 \end{abstract}

\section{Introduction}
\label{section-intro}

    Stochastic gradient descent (SGD) and its variants, e.g., \cite{niu2011hogwild, gorbunov2020unified, lei2020adaptivity, yousefian2012stochastic, nemirovski2009robust, ghadimi2012optimal, mania2017perturbed,  byrd2011use}, are popular and standard methods for large scale optimization and training of various machine learning models, e.g., \cite{bottou2010large, bottou2018optimization, shapiro2021lectures, SchmidtBigData}. Recently, there have been several studies that demonstrate that the gradient noise in SGD that arises, e.g.,
    when training deep learning models, is heavy-tailed, e.g., \cite{simsekli2019heavy, gurbuzbalaban2021heavy, zhang2019adaptive}. 

    Motivated by these studies, we introduce a general analytical framework for \emph{nonlinear} SGD
    when the gradient evaluation is subject to a heavy-tailed noise. We combat the gradient noise with a generic nonlinearity that is applied on the noisy gradient to effectively reduce the noise effect. The resulting class of nonlinear methods subsumes several popular choices in training machine learning models, including normalized gradient descent and clipped gradient descent, e.g., \cite{zhang2019gradient}, the sign gradient, e.g., \cite{balles2020geometry}, and 
    (component-wise) quantized gradient, e.g., \cite{horvath2019stochastic}.\footnote{Interestingly, some of these nonlinear methods are usually introduced with a different motivation than robustness, like, e.g., speeding up training, see, e.g., \cite{zhang2019gradient}, or communication efficiency, \cite{balles2020geometry, bernstein2018signsgd}.}

    We establish for the considered class of methods several results that demonstrate a high degree of robustness to noise
    under very general assumptions on the nonlinearity and on the gradient noise, assuming
    a strongly convex cost with Lipschitz continuous gradient.
    First, for a nonlinearity with bounded outputs (e.g., a sign, normalized, or clipped gradient) and
    the gradient noise that may have infinite moments of order greater than one, assuming
    that the noise probability density function (pdf) is symmetric,
    we show that the nonlinear SGD converges almost surely to the solution, and, moreover, achieves a global~$O(1/t^\zeta)$
    mean squared error (MSE) convergence rate, where we explicitly  quantify the
    degree $\zeta \in (0,1)$. In the same setting, the linear SGD
    generates a sequence with unbounded variances at each iteration~$t$.
    Furthermore, assuming the gradient noise with finite variance,
    we show -- for the unbounded nonlinearities that are lower bounded by a linear function -- 
    almost sure convergence and he~$O(1/t)$ global MSE rate.

    \noindent Next, for the nonlinearities with bounded outputs that can be decoupled component-wise (e.g., a sign or component-wise clipping), we show under the heavy-tail noise a local (asymptotic) $O(1/t)$ rate in the weak convergence sense. More precisely, we show that the sequence generated by the nonlinear SGD is asymptotically normal and explicitly quantify the asymptotic variance. Finally, we illustrate the results on several examples of the nonlinearity and the gradient noise pdf, highlighting and quantifying the noise regimes and the corresponding gains of the nonlinear SGD over the linear SGD scheme. In more detail, the asymptotic variance expression reveals an interesting tradeoff that the nonlinearity makes on the algorithm performance: on the one hand, the nonlinearity suppresses the noise effect to a certain degree, but on the other hand it also reduces the ``useful information flow’’ and hence slows down convergence with respect to the noiseless case. We explicitly quantify this tradeoff and demonstrate through examples that an appropriately chosen nonlinearity strictly improves performance over the linear scheme in a high noise setting. Finally, we carry out numerical experiments on several real data sets that exhibit heavy tail gradient noise effects. The experiments show that, while our analytical framework is more general than usual studies of SGD under heavy-tail noise, several  easy-to-implement example nonlinearities of our framework -- including those not previously used -- are competitive with state of the art alternatives.

    Technically, for component-wise nonlinearities and the asymptotic analysis, we develop proofs based on stochastic approximation arguments, e.g., \cite{nevelson1976stochastic}, following the noise and nonlinearities assumptions framework similar to~\cite{polyak1979adaptive}. The paper \cite{polyak1979adaptive} is concerned with a related but different problem than ours: it considers linear estimation of a vector parameter observed through a sequence of scalar observation equations, and it is not concerned with a global MSE rate analysis that we provide here. For the MSE analysis and for the nonlinearities that cannot be expressed  component-wise, like the clipped and normalized gradient, we develop novel analysis techniques.

    There have been several works that study robustness of stochastic gradient descent under certain variants of heavy-tailed noises. Reference~\cite{zhang2019adaptive} consider an adaptive gradient clipping method and establish convergence rates in expectation for the considered method under a heavy-tailed noise. For this, the authors assume that the expected value of the norm of the gradient noise raised to power $\alpha$ is finite, for $\alpha \in (1,2]$. They also provide lower complexity bounds for SGD methods assuming in addition that the expected $\alpha$-power of the norm of the \emph{stochastic gradient} is finite. The authors of~\cite{gorbunov2020stochastic} consider an accelerated SGD with gradient clipping. They establish high probability bounds for the considered method under the noise that has finite second moment but that does not have to satisfy the sub-Gaussianity assumption. Reference~\cite{davis2021low} proposes a method called proxBoost and establishes for the method high probability bounds, again assuming a finite noise variance and relaxing the sub-Gaussianity assumption. The paper~\cite{simsekli2019heavy} establishes convergence of the \emph{linear} SGD assuming that the gradient noise follows a heavy-tailed $\alpha$-stable distribution. In summary, with respect to existing work, our framework establishes results for the more general setting with respect to both the adopted nonlinearity in SGD and the ``thickness'' of the gradient noise tail, assuming in addition that the noise pdf is a symmetric function. For example, current works usually assume a single choice for the nonlinearity, e.g., gradient clipping, while we consider a general nonlinearity that subsumes many popular choices. Also, provided that the nonlinearity's output is bounded (which is true for many popular choices like the clipped, signed, and normalized gradient), we establish a sublinear MSE convergence rate~$O(1/t^\zeta)$ assuming only that the expected norm of the gradient noise is finite, an assumption weaker than those considered in the works of \cite{gorbunov2020stochastic, zhang2019adaptive, davis2021low, simsekli2019heavy}. On the other hand, we assume a strongly convex smooth cost function, which is equivalent to or stronger than the assumptions made in these works. The MSE rate we establish $\zeta$ may be slower than that in the work of \cite{zhang2019adaptive}. However, the rate $\zeta$ holds uniformly for a general class of nonlinearities, holds for a ``thicker'' noise tail than that in \cite{zhang2019adaptive}, and holds uniformly irrespective of the assumed ``thickness'' of the noise tail. In contrast, the MSE rate in \cite{zhang2019adaptive} holds for a specific adaptive gradient clipping scheme and is dependent on the degree~$\alpha$ of the assumed finite noise moment.

    The idea of employing a nonlinearity into a ``baseline’’ linear scheme has also been used in other contexts. Most notably, several works consider nonlinear versions of the standard consensus algorithm to evaluate average of scalar values in a distributed fashion, e.g., \cite{khan2009distributed, stankovic2019robust, dasarathan2015robust}. The paper~\cite{khan2009distributed} introduces a trigonometric nonlinearity into a standard linear consensus dynamics and shows an improved dependence of the method on initial conditions. References  \cite{stankovic2019robust} and \cite{dasarathan2015robust} employ a general nonlinearity in the linear consensus dynamics and show that it improves the method’s resilience to additive communication noise. The authors of~\cite{sundaram2015consensus} modify the linear consensus by taking out from the averaging operation the maximal and minimal estimates among the estimates from all neighbors of a node. The above works are different from ours as they focus on the specific consensus problem that can be translated into minimizing a convex quadratic cost function in a distributed way over a generic, connected network. In contrast, we consider general strongly convex costs, and we are not directly concerned with distributed systems.

	\textbf{Paper organization}. Section~\ref{section-model-and-framework} describes the problem model and the nonlinear SGD framework that we assume. Section~\ref{section-main-results} and \ref{section-joint-nonlins-iid-noise} explain our results on nonlinear SGD for component-wise and joint nonlinearities, respectively. 
	Sections~\ref{section-proofs-component-wise} and \ref{section-proofs-joint} then provide proofs of the corresponding results. Section~\ref{section-experiments} illustrates the performance of several example methods from our nonlinear SGD framework on real data sets that have heavy-tail gradient noise. Finally, Section~\ref{section-conclusion} concludes the paper. Some auxiliary results and proofs are delegated to the Appendix.

	\textbf{Notation}. We denote by $\mathbb R$ and $\mathbb{R}_+$, respectively, the set of real numbers and real nonnegative numbers, and by ${\mathbb R}^m$ the $m$-dimensional
	Euclidean real coordinate space. We use normal (lower-case or upper-case) letters for scalars,
	lower-case boldface letters for vectors, and upper case boldface letters for
	matrices. Further, we denote by: ${a}_i$ or $[\mathbf{a}]_i$, as appropriate, the $i$-th element of vector $\mathbf{a}$; $\mathbf{A}_{ij}$ or $[\mathbf{A}]_{ij}$, as appropriate, the entry in the $i$-th row and $j$-th column of
	a matrix $\mathbf{A}$;
	$\mathbf{A}^\top$ the transpose of a matrix $\mathbf{A}$; 
	and $\mathrm{trace}(\mathbf{A})$ the sum of diagonal elements of $\mathbf{A}$. Further, we use either  
	$\mathbf{a}^\top \mathbf{b}$ or 
	$\langle \mathbf{a},\,\mathbf{b}\rangle$ 
	for the inner product of vectors 
	$\mathbf{a}$ and $\mathbf{b}$. Next, we let  
	$\mathbf{I}$ and $\mathbf{0}$ be, respectively, the identity matrix and the zero matrix; 
	$\|\cdot\|=\|\cdot\|_2$ the Euclidean (respectively, spectral) norm of its vector (respectively, matrix) argument;  $\phi^\prime(w)$ the first derivative evaluated at $w$ of a function $\phi: {\mathbb R} \rightarrow {\mathbb R}$; $\nabla h(\mathbf{w})$ and $\nabla^2 h(\mathbf{w})$ the gradient and Hessian, respectively, evaluated at $\mathbf{w}$ of a function $h: {\mathbb R}^m \rightarrow {\mathbb R}$; $\mathbb P(\mathcal A)$ and $\mathbb E[u]$ the probability of
	an event $\mathcal A$ and expectation of a random variable $u$, respectively; and by $\textrm{sign}(a)$ the sign function, i.e., $\textrm{sign}(a)=1$, for $a>0$, $\textrm{sign}(a)=-1$, for $a<0$, and $\textrm{sign}(0)=0$.
	Finally, for two positive sequences $\eta_n$ and $\chi_n$, we have: $\eta_n = O(\chi_n)$ if
	$\limsup_{n \rightarrow \infty}\frac{\eta_n}{\chi_n}<\infty$; 
	$\eta_n = \Omega(\chi_n)$ if
	$\liminf_{n \rightarrow \infty}\frac{\eta_n}{\chi_n}>0$; and 
	$\eta_n = \Theta(\chi_n)$ if $\eta_n = O(\chi_n)$ and $\eta_n = \Omega(\chi_n)$.

\section{Problem Model and the nonlinear SGD Framework}
\label{section-model-and-framework}

We consider the following 
 unconstrained problem:
\begin{align}
\label{eq:opt_problem}
\mathrm{minimize\,\,\,\,}f({\mathbf{x}}),
\end{align}
where $f: \mathbb{R}^{d}\mapsto\mathbb{R}$ is a convex
function.

We make the following standard assumption.
\begin{assumption}
\label{assumption-f-i}
Function 
$f: \mathbb{R}^{d}\mapsto\mathbb{R}$ is 
strongly convex with strong convexity parameter 
$\mu>0$, and it has Lipschitz continuous gradient with 
Lipschitz constant~$L \geq \mu$.
\end{assumption}

Under Assumption~\ref{assumption-f-i},
 problem~\eqref{eq:opt_problem} has a unique solution, which we
 denote by $\mathbf{x}^\star \in {\mathbb R}^d$.

In machine learning settings, function $f$ can correspond to
 the risk function, i.e.,
 \begin{equation}
 \label{eqn-risk-fcn}
 f({\mathbf{x}}) = \mathbb{E}_{\,{d} \sim P}\left[\,\ell\left({\mathbf{x}};\mathbf{d}\right) \,\right]+\mathcal{R}({\mathbf{x}}).
 \end{equation}
Here, $P$ is the (unknown) distribution from which the data samples $\mathbf{d} \in {\mathbb R}^q$ are drawn; $\ell(\cdot;\cdot)$
is a loss function, convex in
its first argument
for any fixed value of the second argument;
and $\mathcal{R}:\,{\mathbb R}^d \mapsto \mathbb R$
 is a strongly convex regularizer.
 Similarly,
 $f$ can be empirical risk, i.e.,
 $f({\mathbf{x}}) = \frac{1}{n}\left(\,\sum_{j=1}^{n}
 \ell\left({\mathbf{x}};\mathbf{d}_{j}\right)\,\right) +\mathcal{R}({\mathbf{x}})$,
 where $\mathbf{d}_{j}$,
 $j=1,...,n$, is the set of training data points. 
 Several machine learning models fall within 
 the described framework under Assumption~\ref{assumption-f-i}, 
 including, e.g., $\ell_2$-regularized quadratic and logistic losses.

We introduce a general framework for \emph{nonlinear} SGD methods to solve problem~\ref{eq:opt_problem}; an algorithm within the framework takes the following form:
\begin{equation}
\label{eqn-SGD-nonlinear}
    {\mathbf{x}}^{t+1} = {\mathbf{x}}^t - \alpha_t \boldsymbol{\Psi}(\nabla f({\mathbf{x}}^t) + \boldsymbol{\nu}^t).
\end{equation}
Here, $\mathbf{x}^{t}$ denotes the solution estimate at iteration~$t$, $t=0,1,...$;
$\boldsymbol{\Psi} : \mathbb{R}^d \mapsto \mathbb{R}^d$
is a general nonlinear map; $\alpha_t$ is the employed step size;
$\boldsymbol{\nu}^t \in \mathbb{R}^d$ is a zero-mean gradient noise;
and $\mathbf{x}^0$ is an arbitrary deterministic point in
${\mathbb R}^d$.

 We will specify further ahead the assumptions that we make on the step size $\alpha_t$, the map $\boldsymbol{\Psi}$ and the noise $\boldsymbol{\nu}^t$. 
  Some examples of commonly used maps 
  $\boldsymbol{\Psi}$ that fall within our framework 
  are the following:
 \begin{enumerate}
 \item Sign gradient: $[\boldsymbol{\Psi}(\mathbf{w})]_i = \mathrm{sign}(w_i)$, $i=1,...,d$;
 \item Component-wise clipping: $[\boldsymbol{\Psi}(\mathbf{w})]_i = {w}_i$, for
 $|{w}_i| \leq m$; $[\boldsymbol{\Psi}(\mathbf{w})]_i = m$, for  ${w}_i > m$,
 and  $[\boldsymbol{\Psi}(\mathbf{w})]_i = -m$, for  ${w}_i < -m$, for some constant $m>0$.
 \item Component-wise quantization:
 for each $i=1,...,d$, we let $[\boldsymbol{\Psi}(\mathbf{w})]_i = {r}_j$,
  for $w_i \in (q_{j-1},q_j]$, $j=1,...,J$, where $-\infty =q_0<q_1<...<q_J=+\infty$, $J$ is a positive integer, and the $r_j$'s and
  $q_j$'s are chosen such that each component nonlinearity is an odd function, i.e., $[\boldsymbol{\Psi}(\mathbf{w})]_i =
  -[\boldsymbol{\Psi}(\mathbf{-w})]_i$, for each $i$ and for each $\mathbf{w}$;
 \item Normalized gradient: $\boldsymbol{\Psi}(\mathbf{w}) = \frac{\mathbf{w}}{\|\mathbf{w}\|}$, for $\mathbf{w} \neq 0$, and $\boldsymbol{\Psi}(0)=0$;
 \item Clipped gradient:
 $\boldsymbol{\Psi}(\mathbf{w}) = \mathbf{w}$, for
 $\|\mathbf{w}\| \leq M$, and $\boldsymbol{\Psi}(\mathbf{w}) = \frac{\mathbf{w}}{\|\mathbf{w}\|} \, M$, for  $\|\mathbf{w}\| > M$, for some constant $M>0$.
 \end{enumerate}

Other nonlinearity choices are also introduced ahead (see Section~\ref{section-experiments}).

We next discuss the various possible sources of the gradient noise $\boldsymbol{\nu}^t$. First, the noise may arise due to utilizing a search direction with respect to a data sample. That is, a common search direction in machine learning algorithms is the gradient of the loss with respect to a single data point $\mathbf{d}_i$\footnote{Similar considerations hold for a loss with respect to a mini-batch of data points; this discussion is abstracted for simplicity.}: $\mathbf{g}_i({\mathbf{x}}) = \nabla \ell\left({\mathbf{x}};\mathbf{d}_i\right) + \nabla \mathcal{R}({\mathbf{x}})$. In case of the risk function \eqref{eqn-risk-fcn}, $\boldsymbol{d}_i$ is drawn from distribution~$P$; in case of the empirical risk,  $\boldsymbol{d}_i$ can be, e.g., drawn uniformly at random from the set of data points $\boldsymbol{d}_{j}$, $j=1,...,n$, with repetition along iterations. In both cases, the corresponding gradient noise equals ${\boldsymbol{\nu}} = \mathbf{g}_i({\mathbf{x}}) - \nabla f_i({\mathbf{x}})$. Several recent studies indicate that noise  $\boldsymbol{\nu}$ exhibits heavy tails on many real data sets, e.g, \cite{simsekli2019heavy, gurbuzbalaban2021heavy, zhang2019adaptive}, (See also Section~\ref{section-experiments}.)
  
 We also comment on other possible sources of gradient noise. The noise may be added on purpose to
 the gradient $\nabla f(\mathbf{x})$ for
 improving privacy of an SGD-based learning process, e.g., \cite{pichapati2019adaclip}. 
 Also, the noise $\boldsymbol{\nu}^t$
 may model random computational perturbations or inexact calculations in evaluating
 a gradient~$\nabla f(\mathbf{x})$.

\section{Main results: Component-wise Nonlinearities}
\label{section-main-results}
Section~\ref{section-main-results} provides analysis of the nonlinear SGD method 
for component-wise nonlinearities. 
 That is, we consider here maps $\boldsymbol{\Psi} : \mathbb{R}^d \mapsto \mathbb{R}^d$ of the form 
$\boldsymbol{\Psi}(w_1,..., w_d))^\top$ $= (\Psi(w_1),..., \Psi(w_d))^\top$,
for any $\mathbf{w} \in {\mathbb R}^d$, where (somewhat abusing notation) 
we denote by $\Psi: \,{\mathbb R} \mapsto \mathbb R$
 the component-wise nonlinearity.
 In this setting, we establish for~\eqref{eqn-SGD-nonlinear} almost sure convergence and evaluate
the MSE convergence rate and the asymptotic covariance of the method. 

In more detail, we consider algorithm~\eqref{eqn-SGD-nonlinear}  under Assumptions~\ref{assumption-nu-comp-wise-iid-noise} and~\ref{assumption-psi-comp-wise-iid-noise} below; they follow  
the noise and nonlinearity framework similar to~\cite{polyak1979adaptive}.
\begin{assumption}[Gradient noise] 
\label{assumption-nu-comp-wise-iid-noise}
For the gradient noise random vector sequence $\{\boldsymbol{\nu}^t\}$ in~\eqref{eqn-SGD-nonlinear}, 
$t=0,1,...$,  
$\boldsymbol{\nu}^t \in \mathbb{R}^d$, we assume the following:
\begin{enumerate}
    \item $\{\boldsymbol{\nu}^t\}$ is 
    independent identically distributed (i.i.d.) across iterations, and, for any fixed $t$, $\boldsymbol{\nu}^t$ is independent of ${\mathbf{x}}^t$. Also, random variables $\boldsymbol{\nu}_i^t$ are mutually independent across $i=1,...d$;
    \item Each component ${\nu}_i^t$, $i=1,...,d$, of vector $\boldsymbol{\nu}^t = (\nu_1^t,..., \nu_d^t)^\top$ 
    has a probability density function $p(u)$, $p:\mathbb{R} \mapsto \mathbb{R}_{+}$. The pdf $p$ is symmetric, i.e., $p(u) = p(-u),$ for any $u \in \mathbb{R}$ with $\int |u|p(u) du <+\infty$. 
    \item The pdf $p(u)$ is strictly unimodal, i.e., $p(0) < + \infty$ and $p(v_1) < p(v_2)$ for $|v_1| > |v_2|$.
    \item Function $\Psi$ is strictly increasing, and, for the cumulative distribution function (cdf) associated with pdf $p$, $\Phi(u) = \int_{-\infty}^u p(v)dv$, it holds that $\Phi$ and $\Psi$ have a common growth point, i.e., $\Phi(v + \epsilon) > \Phi(v - \epsilon) $ and $\Psi(v + \epsilon) > \Psi(v - \epsilon)$ for a certain $v \in \mathbb{R}$ and for all $\epsilon > 0$.
    \end{enumerate}
    We assume that at least one of  the conditions 3. or 4. hold.
\end{assumption}
Conditions~1. and~2. in Assumption~\ref{assumption-nu-comp-wise-iid-noise} concerning the requirement that the noise vector is i.i.d. across its components $i=1,...,d$ may be restrictive. For the global MSE analysis,  these assumptions can be relaxed; see ahead the remark after Theorem~\ref{theorem-MSE-comp-Wise} and Appendix~C.

\begin{assumption}[Nonlinearity $\Psi$]
\label{assumption-psi-comp-wise-iid-noise}
Function $\Psi : \mathbb{R} \mapsto \mathbb{R}$ has the following properties:
\begin{enumerate}
\item Function $\Psi$ is a continuous (except possibly on a point set with Lebesgue measure of zero), piece-wise differentiable, monotonically nondecreasing and odd function, i.e., $\Psi(-w) = -\Psi(w),$ for any $w \in \mathbb{R}$;
\item $|\Psi(w)| \leq C_1 \,(1+\,|w|)$, for any $w \in \mathbb R$, for 
some constant $C_1 > 0$.
\item $|\Psi(w)| \leq C_2$, 
for some constant $C_2>0$.
\end{enumerate}
Here, we assume that either 2. or 3. holds.
If 2. holds, then we additionally require a finite 
 variance for the gradient noise, i.e., there holds $\int |u|^2 p(d) du < +\infty$.\footnote{As it will be seen in subsequent text, several statements of results and several proofs treat separately the following two scenarios: 1) condition 2. in Assumption~3 holds (and the gradient noise may not have finite variance); and 2) condition 3. in Assumption~3 holds, but the gradient noise has finite variance. We clearly indicate ahead when we want to distinguish between the two scenarios. For example, when we say that condition 3. in Assumption~3 holds, we refer to the second scenario above. See, e.g., Theorem~\ref{theorem-MSE-comp-Wise}. 
 If the result holds for either of the two scenarios, we do not make specific mention of condition 2. or 3. in Assumption~3. See, e.g., Theorem~\ref{theorem-almost-sure-component-wise}
.}
\end{assumption}

Note that, provided that condition~3. in Assumption~\ref{assumption-psi-comp-wise-iid-noise}  holds, we require only a finite first moment of the gradient noise, while the moments
of $\alpha$-order, $\alpha>1$, may be infinite, hence allowing for heavy-tail
 noise distributions. For example, the gradient noise variance can be infinite. Condition 3. in  Assumption~\ref{assumption-psi-comp-wise-iid-noise}  holds for
several interesting component-wise nonlinearities, like, e.g.,
the sign gradient, component-wise clipping, and quantization schemes 
introduced in Section~\ref{section-model-and-framework}. Note also that Assumption~\ref{assumption-psi-comp-wise-iid-noise} encompasses a broad range of component-wise nonlinearities, beyond the examples in Section~\ref{section-model-and-framework}. (For example, see Section~\ref{section-experiments} for the $\mathrm{tanh}$ and a bi-level quantization nonlinearity.)

Let us define function $\phi: \mathbb{R} \mapsto \mathbb{R},$
as follows. For a fixed (deterministic)  point $w \in \mathbb R$,
$\phi(w)$ is defined by:
\begin{equation}
\label{eqn-phi-definition}
\phi(w) = \mathbb{E} \left[ \Psi(w+\nu_1^0)\right] = \int \Psi(w+u) p(u) du,
\end{equation}
where the expectation is taken with respect to
the distribution of a single entry of the gradient noise at any iteration, i.e., with respect to pdf~$p(u)$. Intuitively, the nonlinearity $\phi$ is a convolution-like transformation  of the nonlinearity $\Psi$, where the convolution is
taken with respect to the gradient noise pdf $p(u)$.
 As we will see ahead,
 the nonlinearity $\phi$ plays
 an effective role in determinining the performance of algorithm~\eqref{eqn-SGD-nonlinear}.

We have the following Theorem.

\begin{theorem}[Almost sure convergence: Component-wise nonlinearity]
\label{theorem-almost-sure-component-wise}
Consider algorithm~\eqref{eqn-SGD-nonlinear}
for solving optimization problem~\eqref{eq:opt_problem},
and let Assumptions \ref{assumption-f-i}--\ref{assumption-psi-comp-wise-iid-noise} hold.
 Assume in addition that $f$ is twice continuously differentiable. 
Further, let the step-size sequence $\{\alpha_t\}$ be square summable, non-summable:
$\sum \alpha_t = +\infty$;
$\sum \alpha_t^2 < +\infty$.
Then, the sequence of iterates $\{\mathbf{x}^t\}$
 generated by algorithm~\eqref{eqn-SGD-nonlinear}
 converges almost surely to the solution $\mathbf{x}^\star$
  of the optimization problem \eqref{eq:opt_problem}.
\end{theorem}

Theorem \ref{theorem-almost-sure-component-wise} establishes a.s. convergence of the nonlinear SGD scheme~\eqref{eqn-SGD-nonlinear}
 under a general setting for the component-wise nonlinearities and gradient noise. For example, provided that the output of the nonlinearity $\Psi$ is bounded, algorithm~\eqref{eqn-SGD-nonlinear} converges
 even when the gradient noise may not have a finite
 $\alpha$-moment, for any $\alpha>1$. (Hence it may have an infinite variance). In contrast, as shown in 
 Appendix~\ref{section-appendixB}, the linear SGD (algorithm \eqref{eqn-SGD-nonlinear} with $\boldsymbol{\Psi}$ 
 being the identity function) generates a sequence 
 of solution estimates with infinite variances, provided that the variance of $p(u)$ is infinite.

\begin{figure}[thpb]
      \centering
      \includegraphics[width=0.7\textwidth, angle =0]{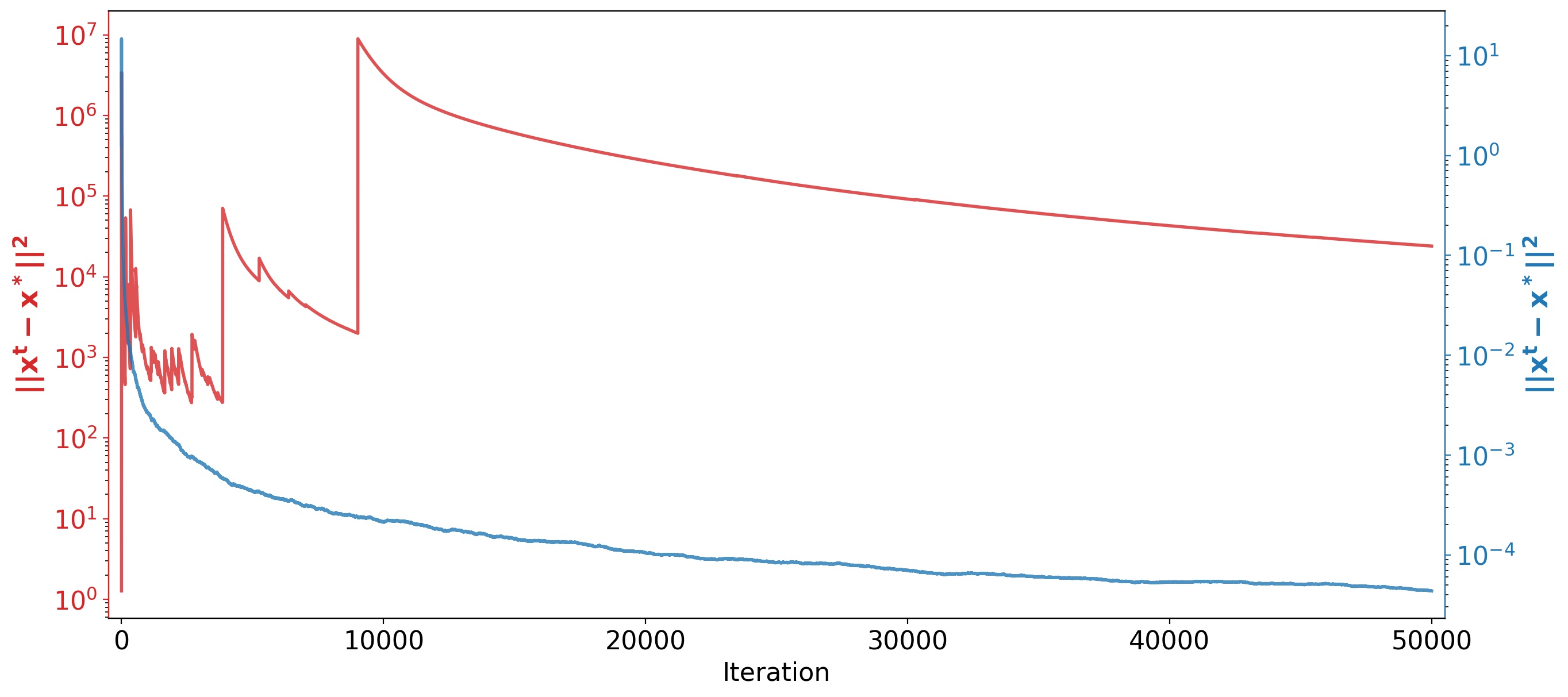}
      \caption{Illustration of Theorem~\ref{theorem-almost-sure-component-wise}: estimated MSE versus iteration counter for the nonlinear SGD in \eqref{eqn-SGD-nonlinear} with component-wise sign nonlinearity (blue line) and the linear SGD (red line).}
      \label{Figure_1}
\end{figure}

\begin{example} Figure \ref{Figure_1} illustrates Theorem~\ref{theorem-almost-sure-component-wise} with a simulation example.  
We consider a strongly convex quadratic function $f:\,{\mathbb R}^d \mapsto \mathbb R$, $f(\mathbf{x}) = \mathbf{x}^\top \mathbf{A} \mathbf{x} + \mathbf{b}^\top \mathbf{x}$, where $\mathbf{A} \in \mathbb{R}^{d \times d}$ is a (symmetric) positive definite matrix, $d = 16$, and quantities $\mathbf{A},\mathbf{b}$ are generated at random. 
 We consider algorithm \eqref{eqn-SGD-nonlinear} 
with the component-wise sign nonlinearity and the linear SGD. 
 The gradient noise 
 is i.i.d. across iterations and across components and has the following pdf:  
     \begin{equation}
        \label{eqn-heavy-tail-noise-example}
        p(u) = \frac{\alpha-1}{2(1+|u|)^{\alpha}},
     \end{equation}
     for $u \in \mathbb R$
      and $\alpha>2$.
      Note that the distribution
      \eqref{eqn-heavy-tail-noise-example}
      does not have a finite $\alpha-1$ moment
       and has finite moments of $r$-th
       order for $r<\alpha-1$. We set in simulation  $\alpha=2.05$. Note that, in this case, the gradient noise has infinite variance.  
  We initialize both the linear and nonlinear algorithm with $\mathbf{x}^0 = 0$, and we let step size $\alpha_t=\frac{1}{t+1}$. 
  Figure~1 shows an estimate of MSE, i.e., of the quantity 
  $\mathbb{E}[\|\mathbf{x}^t-\mathbf{x}^\star\|^2]$,  obtained by averaging results from 100 sample paths. The red line corresponds to the linear SGD, while the blue line corresponds to the nonlinear SGD with the component-wise sign nonlinearity.  As predicted by Theorem~\ref{theorem-almost-sure-component-wise}, the nonlinear SGD drives the MSE to zero, while the linear SGD does not seem to provide a meaningful solution estimate sequence.
\end{example}
\noindent We next establish the mean square error (MSE)
 convergence rate of algorithm~\eqref{eqn-SGD-nonlinear}.

\noindent\begin{theorem}[MSE convergence: Component-wise nonlinearity]
\label{theorem-MSE-comp-Wise}
Consider algorithm~\eqref{eqn-SGD-nonlinear}
for solving optimization problem~\eqref{eq:opt_problem},
and let Assumptions \ref{assumption-f-i},
\ref{assumption-nu-comp-wise-iid-noise}, and  Assumption~\ref{assumption-psi-comp-wise-iid-noise} with condition~3. hold  
Further, let the step-size sequence $\{\alpha_t\}$ be $\alpha_t = a/(t+1)^{\delta}$, $a>0$,
$\delta \in [0.5,1)$.
Then, for the sequence of iterates $\{\mathbf{x}^t\}$
 generated by algorithm~\eqref{eqn-SGD-nonlinear},
 it holds that $\mathbb{E} \left[ \|\mathbf{x}^t-\mathbf{x}^\star\|^2 \right] =O(1/t^\zeta)$,
 or equivalently,
 $\mathbb{E} \left[ f(\mathbf{x}^t) - f^\star \right]  = O(1/t^\zeta)$.
 Here, $\zeta<1$ is any positive number
  such that $\zeta< \mathrm{min}
  \left( 2\delta-1, \frac{a\,(1-\delta)\xi\,\phi^{\prime}(0)\mu}
  {L \,(a\,C_2\,\sqrt{d}\,
 +\|\mathbf{x}^0\|+\|\mathbf{x}^{\star}\|)}\right),$
 and constant $\xi>0$ is such that
 $\phi(a) \geq \frac{\phi^{\prime}(0)\xi}{2}\,a$, for any $a \in (0,\xi)$.
 Furthermore,  let Assumptions \ref{assumption-f-i},
\ref{assumption-nu-comp-wise-iid-noise}, and  Assumption~\ref{assumption-psi-comp-wise-iid-noise} with condition~2. hold, let
$\alpha_t=\frac{a}{(t+1)^{\delta}}$, $\delta\in [0.5,1]$, and assume that
$\inf_{a \neq 0}\frac{|\Psi(a)|}{|a|}>0$.
 Then, there holds that $\mathbb{E} \left[ \|\mathbf{x}^t-\mathbf{x}^\star\|^2 \right] =O(1/t^{\delta})$,
 or equivalently,
 $\mathbb{E} \left[ f(\mathbf{x}^t) - f^\star \right] = O(1/t^{\delta})$. In particular, for $\delta=1$,
 we obtain the $O(1/t)$ MSE rate.
\end{theorem}

\textbf{Remark}. The MSE convergence $O(1/t^{\zeta^\prime})$, for some $\zeta^\prime \in (0,1)$, continues to hold under the same set of assumptions as in Theorem~\ref{theorem-MSE-comp-Wise} but with a relaxed version of  Assumption~\ref{assumption-nu-comp-wise-iid-noise}, where we no longer require that the gradient noise vector has mutually independent components. More precisely, we allow for an i.i.d. noise vector sequence $\{\boldsymbol{\nu}^t\}$, $\boldsymbol{\nu}^t \in {\mathbb R}^d$, that has a symmetric joint pdf $p:\,{\mathbb R}^d \mapsto \mathbb R$, $p(\mathbf{u})=p(-\mathbf{u})$, for any $\mathbf{u} \in {\mathbb R}^d$. In that case, effectively, the role of function $\phi$ in Theorem~\ref{theorem-MSE-comp-Wise} is replaced by functions $w\mapsto \phi_i(w)$, $w \in \mathbb R$, $i=1,...,d$, where $\phi_i(w) = \int \Psi(w+u)p_i(u)du$, and $p_i:\,\mathbb R \mapsto \mathbb R$ is the marginal pdf of the $i$-th component associated with the joint pdf $p:\,{\mathbb R }^d \mapsto \mathbb R$. (See Appendix~C.)

For the bounded nonlinearity case (e.g., sign gradient, component-wise clipping, quantization nonlinearity) and the heavy-tail noise (only the first noise moment assumed to be finite), the nonlinear SGD \eqref{eqn-SGD-nonlinear} achieves a global sublinear MSE rate~$O(1/t^\zeta)$, $\zeta \in (0,1)$. On the other hand, for the finite variance case and an unbounded nonlinearity, the  nonlinear SGD \eqref{eqn-SGD-nonlinear} achieves a global MSE rate $O(1/t)$ provided that $\inf_{w \neq 0}\frac{|\Psi(w)|}{|w|}>0$. This is the best achievable rate and equal to that of the linear SGD in the same setting. Furthermore, by Theorem~\ref{theorem-asymptotic-normality-component-wise} ahead, the nonlinear SGD \eqref{eqn-SGD-nonlinear} with bounded outputs under the heavy-tail noise achieves \emph{locally}, in the weak convergence sense, the faster~$O(1/t)$ rate. This is again in the setting where the linear SGD fails.

\begin{example} We next illustrate 
the value~$\zeta$ in Theorem~\ref{theorem-MSE-comp-Wise} on the family of heavy-tailed pdfs
 in \eqref{eqn-heavy-tail-noise-example}. To be specific, consider the sign nonlinearity $\Psi(w)=\mathrm{sign}(w)$. Then, it is easy to show that:
  \begin{equation*}
  \phi(w) =2 \int_{0}^w p(u) du,~\phi^\prime(0) = 2\,p(0),~\xi \geq  2^{1/\alpha} -1 \approx \frac{1}{\alpha}.
  \end{equation*}
 Using the above calculations, we can see that, for a large $a$, $\zeta$ can be approximated as~$\mathrm{min}\left(2\delta-1,\frac{\mu}{L}\frac{1-\delta}{\sqrt{d}} \right)>0$, $\delta\in (0.5,1)$. 
 
 We also compare the rate~$\zeta$ with the analysis in \cite{zhang2019adaptive} that is closest to our setting with respect to existing work. The authors of \cite{zhang2019adaptive} assume for the MSE upper bound and strongly convex functions analysis that, in our notation, both the quantities $\mathbb{E}[\|\nabla f(\mathbf{x}^t)+\boldsymbol{\nu}^t\|^{\alpha}]$ (a restrictive bounded gradients assumption) and $\mathbb{E}[\|\boldsymbol{\nu}^t\|^{\alpha}]$ are finite for some $\alpha \in (1,2]$. They then show a rate~$O(1/t^{2(\alpha-1)/\alpha})$ for a specific adaptive clipping scheme. This rate can be faster than the~$\zeta$ rate we establish, but ours holds uniformly for a general class of nonlinearities and for any symmetric noise pdf with a finite first moment. 
\end{example}

We next establish asymptotic normality of~\eqref{eqn-SGD-nonlinear}.

\begin{theorem}[Asymptotic normality: Component-wise nonlinearity]
\label{theorem-asymptotic-normality-component-wise}
Consider algorithm~\eqref{eqn-SGD-nonlinear} for solving optimization problem~\eqref{eq:opt_problem}, and let Assumptions \ref{assumption-f-i}--\ref{assumption-psi-comp-wise-iid-noise} hold. Assume in addition that $f$ is twice continuously differentiable. Further, let the step-size sequence $\{\alpha_t\}$
equal: $\alpha_t = a/(t+1)$, $t=0,1,...$, with parameter $a>\frac{1}{2\phi^{\prime}(0)\,\mu}$. Then, the sequence of iterates $\{\mathbf{x}^t\}$ generated by algorithm~\eqref{eqn-SGD-nonlinear} is asymptotically normal, and there holds:
\begin{equation}
    \sqrt{t+1}(\mathbf{x}^t - \mathbf{x}^{\star}) \xrightarrow{d} \mathbb{N}(0,\mathcal{S}),
\end{equation}
where $\xrightarrow{d}$ designates convergence in distribution. The asymptotic covariance $\mathcal{S}$ of the multivariate normal distribution $\mathbb{N}(0,\mathcal{S})$ is given by:
\begin{equation}
\label{eqn-asum-covar-expression}
    \mathcal{S} = a^2\int_{\nu = 0}^{\infty} e^{\nu\boldsymbol{\Sigma}} \mathcal{S}_0e^{\nu\boldsymbol{\Sigma}}d\nu \nonumber
    = a^2\sigma_{\mathrm{\psi}}^2\left[
    2a\phi^{\prime}(0)\nabla^2 f(x^\star) - \mathbf{I}
    \right]^{-1},
\end{equation}
where:
\begin{equation}
    \mathcal{S}_0 = \sigma_{\Psi}^2\,\mathbf{I},\,\,\,\sigma_{\Psi}^2 = \int|\Psi(v)|^2p(v)dv,\,\, 
    \Sigma = \frac{1}{2}\mathbf{I} - a\,\phi^{\prime}(a)\nabla^2f(\mathbf{x}^\star).
\end{equation}
\end{theorem}

Theorem~\ref{theorem-asymptotic-normality-component-wise} establishes asymptotic normality of~\eqref{eqn-SGD-nonlinear} and, moreover, it gives an exact expression for the asymptotic covariance~$\mathcal{S}$ in~\eqref{eqn-asum-covar-expression}, that basically corresponds to the constant in the $1/t$ variance decay near the solution. The asymptotic covariance value \eqref{eqn-asum-covar-expression} reveals an interesting tradeoff with respect to the effect of the nonlinearity $\Psi$. We provide some insights into the tradeoff through examples below.

\begin{example} We compare the linear SGD and the nonlinear SGD with component wise clipping. For illustration and simplification of calculations, we consider the special case when $\nabla^2 f(\mathbf{x}^\star)$ is a symmetric matrix with all eigenvalues equal to one. Then, it is straightforward to show that the per-entry asymptotic variance for the best choice of parameter $a$ over the admissible set of values equals:
 \begin{equation}
  \label{eqn-asymptotic-covariance-closed-form}
  \inf_{a>\frac{1}{2\phi^\prime(0)}} \mathrm{trace}
  \left( \mathcal{S}\right) = \frac{\sigma_{\mathrm{\Psi}}^2}
  {(\phi^{\prime}(0))^2}.
  \end{equation}
Here, for the linear SGD i.e., when $\Psi(a)=a$, we have that $\sigma_{\mathrm{\Psi}}^2 = \int a^2 p(a)da$ equals the gradient noise (per component) variance~$\sigma_{\nu}^2$, and $\phi^{\prime}(0)=1$, and so \eqref{eqn-asymptotic-covariance-closed-form} equals~$\sigma_{\nu}^2$. Now, consider the coordinate-wise clipping, with $\Psi(a)=a$ for $|a| \leq m$ and $\Psi(a)=\mathrm{sign}(a)\,m$, for $|a|>m$, for some $m>0$. Then, we have: $\sigma_{\Psi}^2 = m^2 - 2\,\int_0^m (m^2-v^2)p(v)dv$, and $\phi^\prime(0) = 2\,\int_0^m p(v) dv.$ Note that the case~$m \rightarrow \infty$ corresponds to the linear SGD case. Consider now the tradeoff with respect to the choice of~$m$. Clearly, taking a smaller $m$ has a positive effect on the numerator in \eqref{eqn-asymptotic-covariance-closed-form} (it suppresses the noise effect). On the other hand, reducing $m$ has a negative effect on the denominator in \eqref{eqn-asymptotic-covariance-closed-form}; that is, it reduces the value~$\phi^\prime(0)$ -- intuitively, it  ``lowers the quality'' of the  search direction utilized with \eqref{eqn-SGD-nonlinear}. One needs to choose the nonlinearity, i.e., the parameter~$m$, optimally, to strike the best balance here. Clearly, for larger gradient noise~$\sigma_{\nu}^2$, we should pick a smaller value of~$m$. It can be shown that, for any finite $\sigma_{\nu}^2$, there is an optimal value $m^\star \in (0,\infty)$ that minimizes~\eqref{eqn-asymptotic-covariance-closed-form}.   
\end{example}

\begin{example}
    We continue to assume the simplified setting when the per-entry asymptotic variance equals~\eqref{eqn-asymptotic-covariance-closed-form}.  We consider the sign gradient nonlinearity and the class of heavy-tail gradient noise distributions in \eqref{eqn-heavy-tail-noise-example}. It can be shown that here: $\sigma_{\Psi}^2=1$; $\sigma_{\nu}^2=\frac{2}{(\alpha-3)(\alpha-2)}$, for $\alpha>3$ and $\sigma_{\nu}^2=\infty$, else; and $\phi^\prime(0)=\alpha-1.$ Therefore, for the sign gradient, the best achievable per entry asymptotic variance equals $\frac{1}{(\alpha-1)^2}$, while for the linear SGD it equals $\frac{2}{(\alpha-2)(\alpha-3)}$ for $\alpha>3$, and is infinite for $\alpha \in (2,3]$. Hence, we can see for the considered example that the sign gradient outperforms the linear SGD for any $\alpha>2$, and the gap becomes larger as $\alpha$ gets smaller.
\end{example}

\begin{figure}[thpb]
      \centering
      \includegraphics[width=0.7\textwidth]{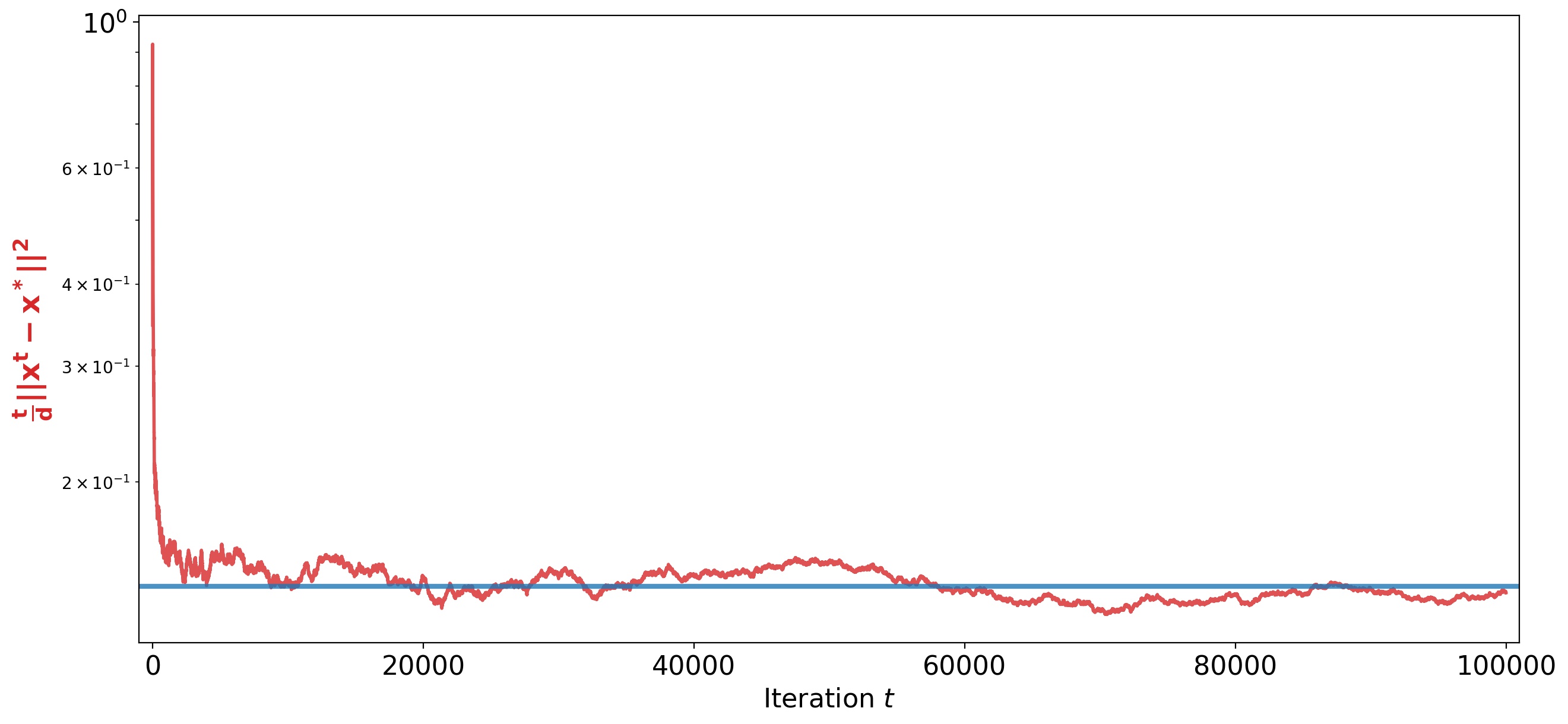}
      \caption{Illustration of Theorem~\ref{theorem-asymptotic-normality-component-wise}: Monte Carlo estimate of per-entry asymptotic variance (red line) and the theoretical per-entry asymptotic variance in Theorem~\ref{theorem-asymptotic-normality-component-wise} (blue line).}
      \label{Figure_2}
\end{figure}

\begin{example}
    We still consider the simplified setting of~\eqref{eqn-asymptotic-covariance-closed-form}. If the noise pdf $p(u)$ is known, then, following~\cite{polyak1979adaptive}, we can find a globally optimal nonlinarity that minimizes~\ref{eqn-mathcal-S} that takes the form: $\Psi(a) = - \frac{d}{da} \mathrm{ln}(p(a))$. The corresponding optimal asymptotic variance equals the Fisher information associated with the pdf~$p(u)$.
\end{example}
\begin{example}
    Figure~\ref{Figure_2} illustrates Theorem~\ref{theorem-asymptotic-normality-component-wise} for the nonlinear SGD in~\eqref{eqn-SGD-nonlinear} with component-wise sign nonlinearity and the same simulation setting used for the numerical illustration of Theorem~\ref{theorem-almost-sure-component-wise} and step-size $\alpha_t = \frac{10}{t+1}$. The red line plots quantity $\frac{t}{{d}}{\|\mathbf{x}^t - \mathbf{x}^\star\|^2}$ estimated through 100 sample path runs. This quantity estimates the constant in the $1/t$ per-entry asymptotic variance decay, i.e., it is a numerical estimate of the per-entry asymptotic variance $\frac{\mathrm{trace}(\mathcal{S})}{d}$, where $\mathcal{S}$ is given in  Theorem~\ref{theorem-asymptotic-normality-component-wise}. The blue horizontal line marks the value  $\frac{\mathrm{trace}(\mathcal{S})}{d}$. We can see that the simulation matches well the theory.  
\end{example}

\section{Main results: Joint Nonlinearities} 
\label{section-joint-nonlins-iid-noise}
We now consider algorithm~\eqref{eqn-SGD-nonlinear} for a nonlinearity $\boldsymbol{\Psi}:\,{\mathbb R}^d\mapsto {\mathbb R}^d$ that cannot be decoupled into (equal) component wise nonlinearities $\Psi: \,{\mathbb R} \mapsto {\mathbb R}$, as it was possible before. More precisely, we make the following assumptions on the gradient noise ${\boldsymbol{\nu}}^t$ and the nonlinear map $\boldsymbol{\Psi}: \,{\mathbb R}^d\mapsto {\mathbb R}^d$.

 \begin{assumption}
\label{assumption-p-joint-iid-noise}[Gradient noise]
 For the gradient noise sequence $\{\boldsymbol{\nu}^t\}$, we assume the following:
\begin{enumerate}
    \item  The sequence of random vectors $\{\boldsymbol{\nu}^t\}$ is i.i.d. and for any $t=0,1,...$, $\boldsymbol{\nu}^t$ is independent of ${\mathbf{x}}^t$. Moreover, $\boldsymbol{\nu}^t$ has a joint symmetric pdf $p(\mathbf{u})$, $p:\,{\mathbb R}^d \mapsto \mathbb R$, i.e., $p(\mathbf{u}) = p(-\mathbf{u})$, for any $\mathbf{u} \in \mathbb{R}^d$ with $\int\|{\mathbf{u}}\|p(\mathbf{u})d\mathbf{u} < \infty$;
     \item 
     There exists a positive constant $B_0$ such that, for any $ {\mathbf{x}} \in {\mathbb R}^d$, ${\mathbf{x}} \neq 0$,
     for any $A \in (0,1],$ there exists $ {\lambda} ={\lambda}(A) > 0$, such that\footnote{The integration set $\{\mathbf{u} \in \mathbb{R}^d:\,\frac{\mathbf{u}^\top \mathbf{x}}{\|\mathbf{u}\|\|\mathbf{x}\|} \in [0,A],\,\|u\| \leq B_0 \}$ also includes the point $\mathbf{u}=0$. In other words, for compact notation here and throughout the paper, we write 
     $\frac{\mathbf{u}^\top \mathbf{x}}{\|\mathbf{u}\|\|\mathbf{x}\|} \in [0,A]$ instead of 
     $0 \leq 
     \mathbf{u}^\top \mathbf{x}
     \leq
     A\,\|\mathbf{u}\|\|\mathbf{x}\|$.}
    $\int_{\{\mathbf{u} \in \mathbb{R}^d:\,\frac{\mathbf{u}^\top \mathbf{x}}{\|\mathbf{u}\|\|\mathbf{x}\|} \in [0,A],\,\|u\| \leq B_0 \}} $ $p(\mathbf{u})d\mathbf{u} > {\lambda}(A).$
\end{enumerate}
\end{assumption}

Assumption~\ref{assumption-p-joint-iid-noise} allows for a heavy-tailed noise 
vector whose components can be mutually dependent. Condition 2. in Assumption~\ref{assumption-p-joint-iid-noise} is mild; it says that the joint pdf $p(\mathbf{u})$ is ``non-degenerate'' in the sense that, along each ``direction'' (determined by arbitrary nonzero vector~$\mathbf{x}$), the intersection of the set $\{\frac{\mathbf{u}^\top \mathbf{x}}{\|\mathbf{u}\|\|\mathbf{x}\|} \in [0,A]\}$ and the ball $\{\|\mathbf{u}\| \leq B_0\}$ consumes a positive mass of the joint pdf $p(\mathbf{u})$.

\begin{assumption}[Nonlinearity $\boldsymbol{\Psi}$]
\label{assumption-psi-joint-iid-noise}
The nonlinear map $\boldsymbol{\Psi}: \,{\mathbb R}^d \mapsto {\mathbb R}^d$ takes the folloing form: $\boldsymbol{\Psi}(\mathbf{w}) = \mathbf{w}\mathcal{N}(\|{\mathbf{w}}\|)$, where function $\mathcal{N}:\mathbb{R}_+ \mapsto \mathbb{R}_+$ satisfies the following: 
\begin{enumerate}
  \item Function $\mathcal{N}$ is non-increasing and continuous except possibly on a point set with Lebesgue measure of zero with $\mathcal{N}(q) > 0$, for any  $q > 0$. The function $q \mathcal{N}(q)$ is non-decreasing;
  \item $\|\boldsymbol{\Psi}(\mathbf{w})\| \leq C_1^\prime (1+\|\mathbf{w}\|),$ for any $\mathbf{w} \in \mathbb{R}^d$, for some $C_1^\prime > 0$;
  \item $\|\boldsymbol{\Psi}(\mathbf{w})\| \leq C_2^\prime,$ for any $\mathbf{w} \in \mathbb{R}^d$, for some $C_2^\prime > 0$.
\end{enumerate}

Here, we assume that either 2 or 3 holds. If 2 holds, then we additionally require that the second moment of $\boldsymbol{\nu}^t$ is bounded, i.e., $\int\|\mathbf{u}\|^2 p(\mathbf{u}) d\mathbf{u}<\infty$.\footnote{Analogously to Assumption~3, we often refer to the two different scenarios, corresponding to conditions 2. and 3. in Assumption, 5, respecticvely. See also the footnote at the end of Assumption~3.}
\end{assumption}

There are many nonlinearities that satisfy Assumption \ref{assumption-psi-joint-iid-noise}, including, the normalized gradient and the clipped gradient discussed in Section~\ref{section-model-and-framework}.

\vspace{-0mm}
\begin{figure*}[t!]
\centering
\subfloat[\label{fig:2a}]{\includegraphics[width=0.33\linewidth]{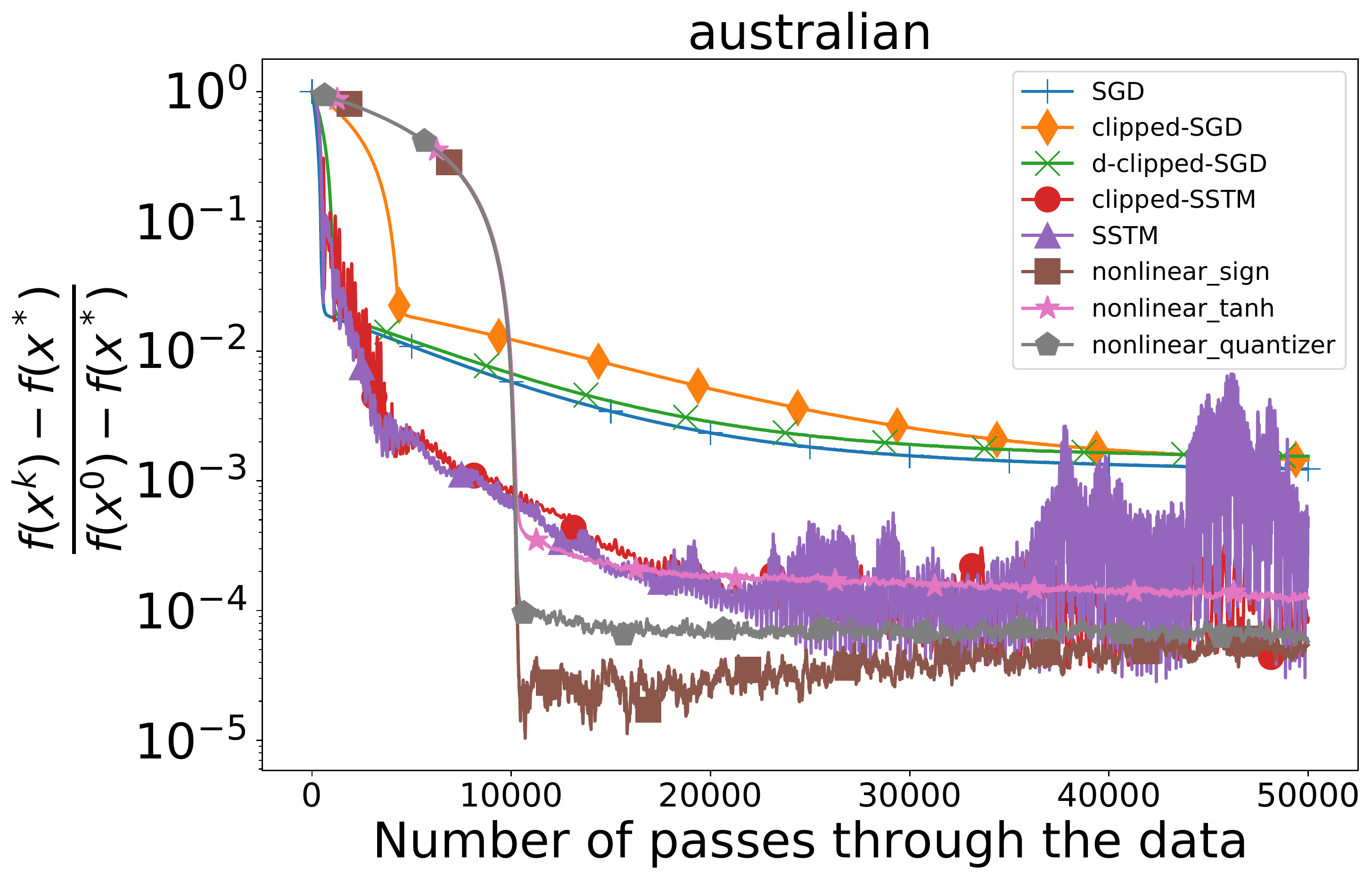}}
\centering
\subfloat[\label{fig:2b}]{\includegraphics[width=0.33\linewidth]{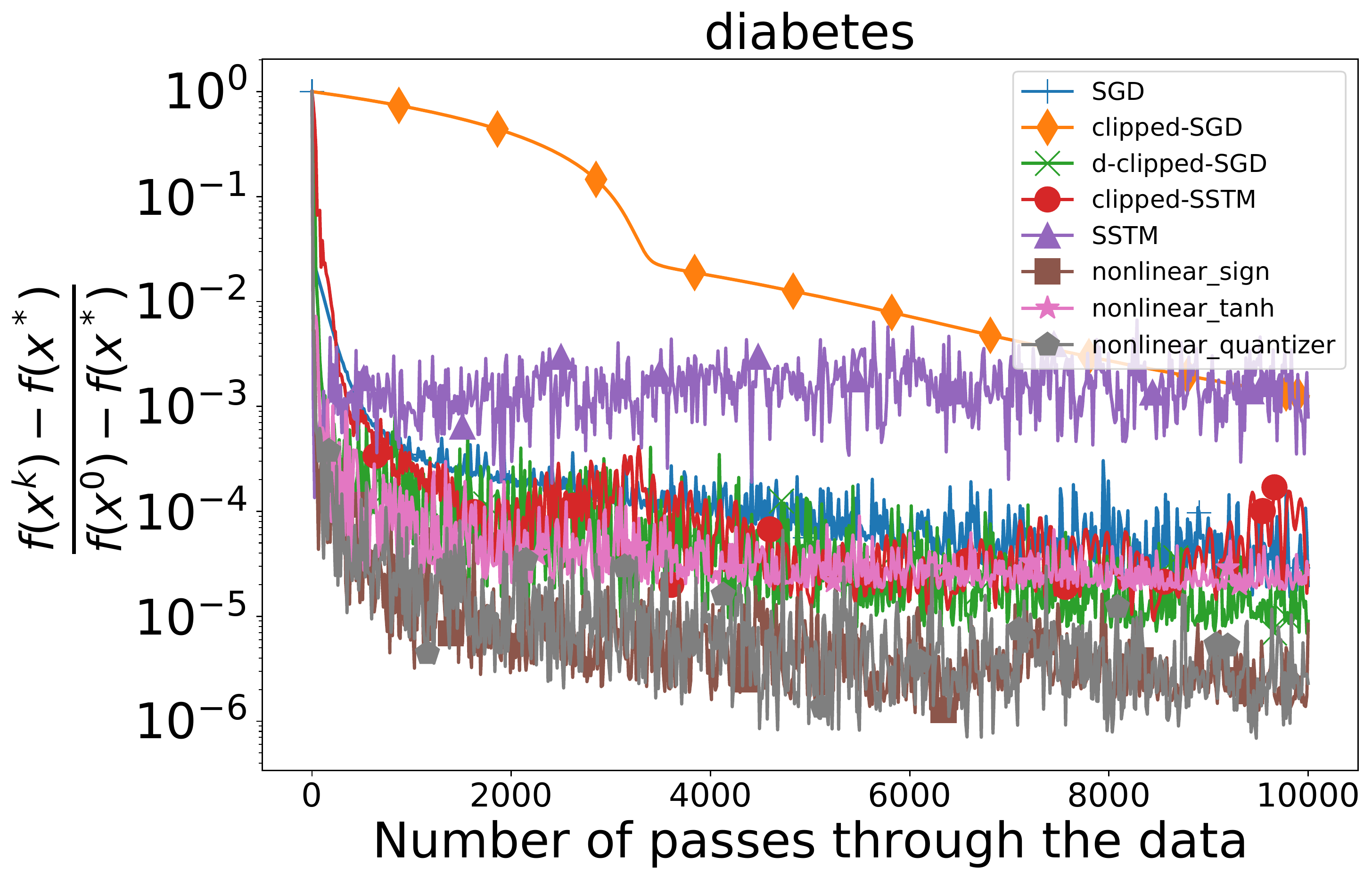}}
\centering
\subfloat[\label{fig:2c}]{\includegraphics[width=0.33\linewidth]{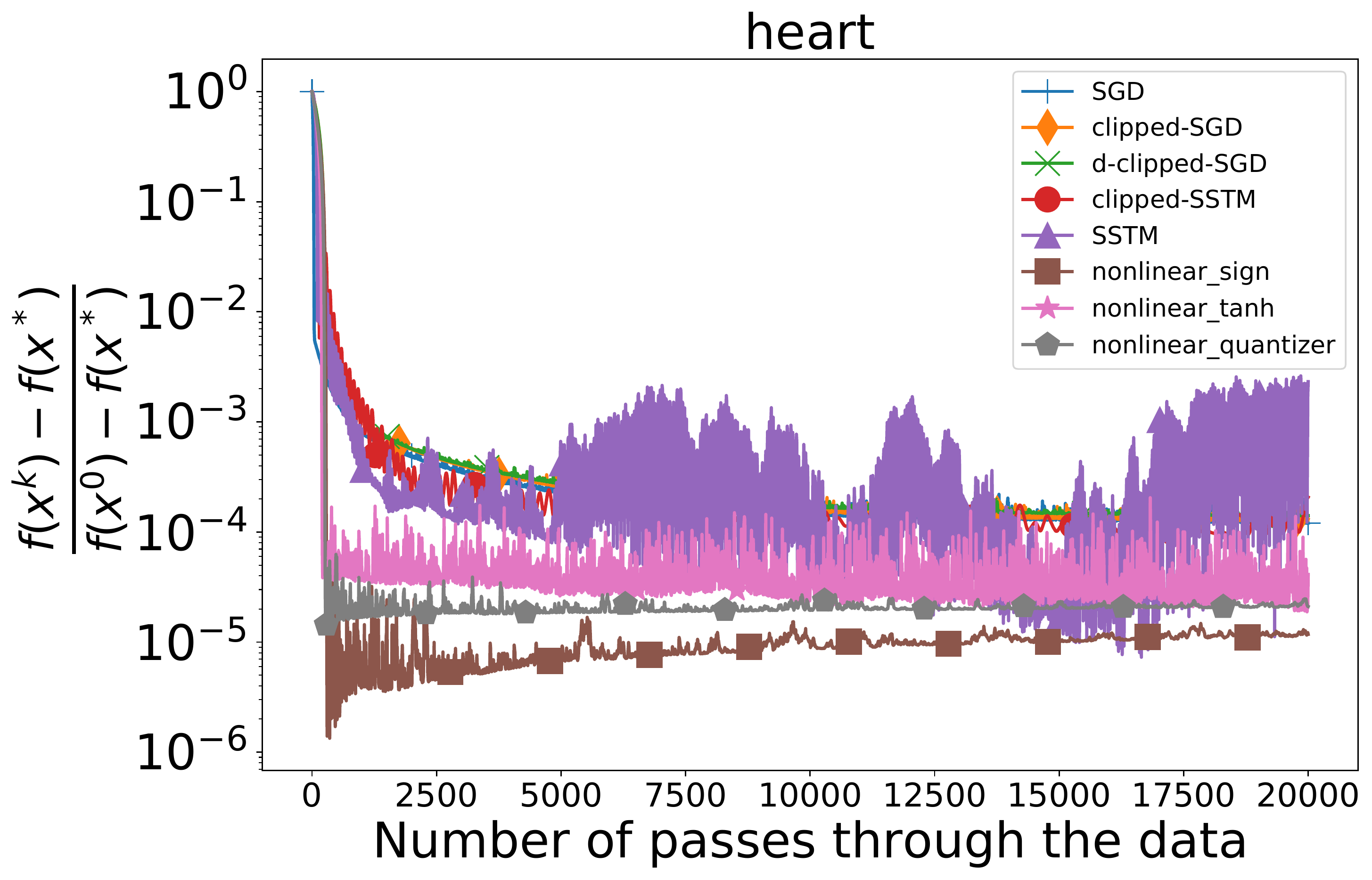}}
\caption{Comparison of the optimization algorithms across different datasets}
\label{fig:heavy_tail_comp}
\end{figure*}
\vspace{-0mm}
\begin{theorem}[MSE and a.s. convergence: Joint nonlinearity]
\label{theorem-MSE-joint-iid}

Consider algorithm~\eqref{eqn-SGD-nonlinear} for solving optimization problem~\eqref{eq:opt_problem}, and let Assumptions \ref{assumption-f-i},
\ref{assumption-p-joint-iid-noise}, and Assumption~\ref{assumption-psi-joint-iid-noise} with condition~3. hold. Further, let the step-size sequence $\{\alpha_t\}$ be $\alpha_t = a/(t+1)^{\delta}$, $a>0$, $\delta \in (0.5,1)$. Then, for the sequence of iterates $\{\mathbf{x}^t\}$ generated by algorithm~\eqref{eqn-SGD-nonlinear}, it holds that $\mathbb{E} \left[ \|\mathbf{x}^t-\mathbf{x}^\star\|^2 \right] =O(1/t^\zeta)$, or equivalently, $\mathbb{E} \left[ f(\mathbf{x}^t) - f^\star \right] = O(1/t^\zeta)$, where $\zeta \in (0,1)$. In alternative, let Assumptions \ref{assumption-f-i}, \ref{assumption-p-joint-iid-noise}, and  Assumption~\ref{assumption-psi-joint-iid-noise} with condition~2. hold, let $\alpha_t=\frac{a}{(t+1)^{\delta}}$, $\delta\in (0.5,1]$, and assume that $\inf_{\mathbf{w} \neq 0}\frac{\|\boldsymbol{\Psi}(\mathbf{w})\|}{\|\mathbf{w}\|}>0$. Then, $\mathbb{E} \left[ \|\mathbf{x}^t-\mathbf{x}^\star\|^2 \right] =O(1/t^{\delta})$, or equivalently, $\mathbb{E} \left[ f(\mathbf{x}^t) - f^\star \right] = O(1/t^{\delta})$. In particular, for $\delta=1$ and a sufficiently large parameter $a$, we obtain the $O(1/t)$ MSE rate. Finally, under Assumptions \ref{assumption-f-i}, \ref{assumption-p-joint-iid-noise}, and~\ref{assumption-psi-joint-iid-noise}, $\delta\in (0.5,1]$ and $f$ that is in addition twice continuously differentiable, we have that $\mathbf{x}^t$ converges to $\mathbf{x}^\star$, a.s.
\end{theorem}

Theorem~\ref{theorem-MSE-joint-iid} establishes a.s. convergence and a global sublinear MSE rate of algorithm~\eqref{eqn-SGD-nonlinear} for a nonlinearity with bounded outputs (e.g., a normalized or clipped gradient) in the presence of heavy-tailed noise that may have infinite moments of order greater than one. See also the discussion after Theorem~\ref{theorem-MSE-comp-Wise} for analogous interpretations and comparisons with existing work. See ahead \eqref{eqn-zeta-joint} in the proof of Theorem~\ref{theorem-MSE-joint-iid} for the obtained bound on rate~$\zeta$.

\section{Intermediate results and proofs: Component-wise nonlinearities}
\label{section-proofs-component-wise}
This section provides proofs of Theorems \ref{theorem-almost-sure-component-wise}--\ref{theorem-asymptotic-normality-component-wise}, accompanied with the required intermediate results. Subsection~5.1 deals with the asymptotic analysis (Theorems \ref{theorem-almost-sure-component-wise} and \ref{theorem-asymptotic-normality-component-wise}), while Subsection 5.2 considers MSE analysis (Theorem \ref{theorem-MSE-comp-Wise}).

\subsection{Asymptotic analysis: Proofs of Theorems \ref{theorem-almost-sure-component-wise} and \ref{theorem-asymptotic-normality-component-wise}}
The next Lemma, due to \cite{polyak1979adaptive}, establishes structural properties of function~$\phi$ in \eqref{eqn-phi-definition}. The Lemma says that essentially, the convolution-like transofrmation of the nonlinearity preserves the structural properties of the nonlinearity.

\begin{lemma} \cite{polyak1979adaptive}
\label{Lemma-Polyak}
Consider function $\phi$ in \eqref{eqn-phi-definition}, 
where function $\Psi:\,\mathbb R \mapsto \mathbb R$ 
satisfies Assumption~3. Then, the following holds.
\begin{enumerate}
    \item $\phi$ is odd;
    \item If $|\Psi(\nu)| \leq C_2,$ for any $\nu \in \mathbb R$, then $|\phi(a)| \leq K_2,$ for any $a \in \mathbb R$, for some constant $K_2>0$;
    \item If $|\Psi(\nu)| \leq C_1 (1 + |\nu|)$, for any $\nu \in \mathbb R$, then $|\phi(a)| \leq K_1 (1 + |a|)$, 
    for any $a \in \mathbb R$, for some constant $K_1>0$;
    \item $\phi(a)$ is monotonically nondecreasing;
    \item $\phi(a) > 0,$ for any $ a>0$.
    \item $\phi$ is continuous at zero;
    \item $\phi$ is differentiable at zero, with a strictly positive derivative at zero, equal to:
    \begin{equation}
    \begin{split}
        \phi^\prime(0) &= \sum_{i=1}^s \left(\Psi(\nu_i + 0) + \Psi(\nu_i - 0)\right)p(\nu_i) + \sum_{i=0}^s\int_{\nu_i}^{\nu_{i+1}}\Psi'(\nu)p(\nu)d\nu ,
    \end{split}
    \end{equation}
    where $\nu_i, i=1,...,s$ are points of discontinuity of $\Psi$ such that $\nu_0 = -\infty$ and $\nu_{s+1} = +\infty$
\end{enumerate}
\end{lemma}

We proceed by setting up the proof of  Theorem~\ref{theorem-almost-sure-component-wise}. The proof relies on convergence analysis of single-time scale stochastic approximation methods from \cite{nevelson1976stochastic}; more precisely, we utilize Theorem~\ref{theorem-29} in the Appendix; see also \cite{kar2012distributed}.

We first put algorithm~\eqref{eqn-SGD-nonlinear} in the format that complies with Theorem~\ref{theorem-29}. Namely, algorithm~\eqref{eqn-SGD-nonlinear} can be written as:

\begin{equation}
\label{eqn-proof-SA-format}
    \mathbf{x}^{t+1} = \mathbf{x}^t + \alpha_t\left[{\mathbf{r}}(\mathbf{x}^t) + {\boldsymbol{\gamma}}(t+1, x^t, \omega)\right].
\end{equation}
Here, $\omega$ denotes an element of the underlying probability space, and
\begin{equation}
\label{eqn-R}
    {\mathbf{r}}(\mathbf{x}) = -\boldsymbol{\phi}(\nabla f(\mathbf{x})),
\end{equation}
where, abusing notation, $\boldsymbol{\phi}: \mathbb{R}^d \mapsto \mathbb{R}^d$ is defined by
 $ \left(\boldsymbol{\phi}(a_1,...,a_d)\right)^\top =$ $\left(\phi(a_1),...,\phi(a_d)\right)^\top$. That is, we have that:
\begin{equation}
    {\mathbf{r}}(\mathbf{x}) = - \left(\phi[\nabla f(x))_1],..., \phi[\nabla f(x))_d] \right)^\top
\end{equation}
and
\begin{equation}
\label{eqn-Gamma}
    {\boldsymbol{\gamma}}(t+1, \mathbf{x}, \omega) = \boldsymbol{\phi}(\nabla f(\mathbf{x})) - \Psi(\nabla f(\mathbf{x}) + \boldsymbol{\nu}^t).
\end{equation}

We provide an intuition behind the algorithmic format~\eqref{eqn-proof-SA-format}. Quantity ${\mathbf{r}}(x)$ is a deterministic, ``useful'', progress direction with respect to the evolution of $\mathbf{x}^t$; quantity ${\boldsymbol{\gamma}}(t+1, x, \omega)$ is the stochastic component that plays a role of a noise in the system.

We adopt the following Lyapunov function: $V(x) = f(x) - f^\star,$ $V:\mathbb{R}^d \mapsto \mathbb{R}$, where $f^\star = \inf_{x\in \mathbb{R}^d}f(x) = f(x^\star)$. We are ready to prove Theorem~\ref{theorem-almost-sure-component-wise}.

\emph{Proof} (Proof of Theorem~\ref{theorem-almost-sure-component-wise}). 
 We now verify conditions B1-B5 from Theorem~\ref{theorem-29}. It can be shown that, under Assumptions~2 and~3, $\phi(a)$ continuous for any $a \in \mathbb R $(see the proof of Lemma~5 and Theorem~1 in \cite{polyak1979adaptive}), and therefore, in view of Assumptions 1--3, for each~$t$, function $\boldsymbol{\gamma}(t+1,\cdot,\cdot)$ is measurable. Hence, condition~B1 holds. Consider the filtration~$\mathcal{H}_t$, $t=1,2,...$, where $\mathcal{H}_t$ is the $\sigma$-algebra generated with random vectors $\boldsymbol{\nu}^s$, $s=0,...,t-1$. Then, the family of random vectors $\left\{{\boldsymbol{\gamma}}\left(t,\mathbf{x},\omega\right)
\right\}_{\mathbf{x}\in\mathbb{R}^{d}}$ is $\mathcal{H}_{t}$ measurable, zero-mean and independent of $\mathcal{H}_{t-1}$. Thus, condition B2 holds.

For B3, we need to prove that
\begin{equation}
\label{eqn-condition-to-prove}
\sup_{\mathbf{x}:\|\mathbf{x}-\mathbf{x}^\star \|\in(\epsilon, \frac{1}{\epsilon})}
\langle {\mathbf{r}}(\mathbf{x}), \frac{\partial V}{\partial x}(\mathbf{x}) \rangle < 0,\, \mathrm{for\,any}\, \epsilon > 0,
\end{equation}
where $\frac{\partial V}{\partial x}(\mathbf{x}) = \nabla f(\mathbf{x}).$ Let us fix an $\epsilon > 0$. Consider arbitrary  $\mathbf{x}\in\{\mathbf{y}\in \mathbb{R}^d:\,\|\mathbf{y}-\mathbf{x}^\star\| \in (\epsilon, \frac{1}{\epsilon})\}.$ Then, we have: 
\begin{eqnarray*}
\langle {\mathbf{r}}(\mathbf{x}), \frac{\partial V}{\partial x} \rangle 
&=& -{\boldsymbol{\phi}}(\nabla f(\mathbf{x}))^\top (\nabla f(\mathbf{x})) \\
&=& -\sum_{j=1}^{d}\phi([\nabla f(\mathbf{x})]_j) [\nabla f(\mathbf{x})]_j  
= -\sum_{j=1}^{d}|\phi([\nabla f(x)]_j)| \,|[\nabla f(x)]_j|,
\end{eqnarray*}

\noindent 
where the last inequality holds because $\phi$ is an odd function. Since $\|\mathbf{x}-\mathbf{x}^\star\| > \epsilon$ and $\|\nabla f(\mathbf{x})\|^2 > \frac{\mu}{2}\|\mathbf{x}-\mathbf{x}^{\star}\|^2$ (due to strong convexity of $f$), we have $\|\nabla f(\mathbf{x})\| > \sqrt{\frac{\mu}{2}}\epsilon$, where we recall that $\mu$ is the strong convexity constant of~$f$. Therefore, there exists an index $i\in\{1,...,d\}$ such that $|[\nabla f(\mathbf{x})]_i| > \frac{1}{d}\sqrt{\frac{\mu}{2}}\epsilon =: \epsilon^{\prime}$. Next, because $\phi^\prime(0) > 0$, and $\phi$ is continuous at 0 and is nondecreasing (by Lemma~\ref{Lemma-Polyak}), we have that $|\phi(b)| > \delta$ for some $\delta=\delta(\epsilon) > 0$, for all $b \in (\epsilon,1/\epsilon)$. Finally, we have that: $\leq - \epsilon'\delta(\epsilon)$, for any $\mathbf{x}$ such that $\|\mathbf{x}-\mathbf{x}^{\star} \|\in(\epsilon, \frac{1}{\epsilon}),$ and therefore $\sup_{\mathbf{x}:\|\mathbf{x}-\mathbf{x}^{\star} \|\in(\epsilon, \frac{1}{\epsilon})}$ $\langle {\mathbf{r}}(x), \frac{\partial V}{\partial x} \rangle < 0$, hence verifying condition B3.

We next verify condition~B4. Consider quantity ${\mathbf{r}}(\mathbf{x})$ in \eqref{eqn-R}. By Lemma~{5.1} and the fact that $f$ has Lipschitz gradient and is strongly convex (Assumption \ref{assumption-f-i}), it follows that:
\begin{equation}
\| {\mathbf{r}}(\mathbf{x})\|^2 \leq C_1 + C_2 V(\mathbf{x}),
\end{equation}
for some positive constants $C_1$ and $C_2$. 
Also, since
\begin{equation}
\|{\boldsymbol{\gamma}}(\mathbf{x}, t+1, \omega)\|^2 \leq
2\|\phi(\nabla f(\mathbf{x}))\|^2 + 2\|\Psi(\nabla f(\mathbf{x}) + \boldsymbol{\nu}^t)\|^2,
\end{equation}
and it holds that either
1)  $\Psi$ is bounded or 2) $|\Psi(a)| \leq C_2\,(1 + |a|)$ and $\nu_i^t$ has a finite variance, we have:
\begin{equation}
\mathbb{E} \left[ \|{\boldsymbol{\gamma}} (\mathbf{x}, t+1, \omega)\|^2 \right ]  \leq C_3 + C_4 \, V(\mathbf{x}),
\end{equation}
for some positive constants $C_3,C_4$.
Now, we finally have:
\begin{equation}
    \| {\mathbf{r}}(\mathbf{x} )\|^2 + \mathbb{E} \left[ \|{\boldsymbol{\gamma}}(\mathbf{x}, t+1, \omega)\|^2 \right ] \leq C_5 + C_6\, V(\mathbf{x}),
\end{equation}
for some positive constants $C_5, C_6$, 
and hence condition B4 holds. 
Condition B5 holds by the choice of the step size sequence
$\{\alpha_t\}$ in the Theorem statement.  
 Summarizing, all
 conditions
B1-B5 hold true, and hence $\mathbf{x}^t \rightarrow \mathbf{x}^{\star}$, almost surely. $\Box$

We continue by proving Theorem~\ref{theorem-asymptotic-normality-component-wise}.

\textit{Proof} (Proof of Theorem~\ref{theorem-asymptotic-normality-component-wise}). 
We prove the Theorem 
 by verifying conditions C1-C5 in  Theorem~\ref{theorem-29}.  
  To verify condition C1, consider
  ${\mathbf{r}}(\mathbf{x})$ in \eqref{eqn-R}
   and note that, using the mean value theorem, it can be expressed as follows:
\begin{equation}
\begin{split}
    {\mathbf{r}}(\mathbf{x}) &= -\phi(\nabla f(\mathbf{x}) - \nabla f(\mathbf{x}^{\star}))\\
    &= -\phi \left ( \underbrace{\left[ \int_0^1 \nabla^2 f(\mathbf{x}^{\star} +
    t(\mathbf{x} - \mathbf{x}^{\star})dt \right]}_{\mathbf{H}_t} (\mathbf{x}- \mathbf{x}^{\star}) \right ) \\
    &= -\phi \left ( \mathbf{H}_t (\mathbf{x} - \mathbf{x}^*) \right ) = - \phi^{\prime}(0) \nabla^2 f(\mathbf{x}^{\star}) (\mathbf{x} - \mathbf{x}^{\star}) + \delta(\mathbf{x}),
\end{split}
\end{equation}
where $\lim_{\mathbf{x}\rightarrow \mathbf{x}^{\star}}\frac{\|\delta(x)\|}{\|x-x^*\|} = 0$.
 Hence, in the notation of Theorem~\ref{theorem-29},
 we have that $\mathbf{B} = -\phi^\prime (0) \nabla^2 f(\mathbf{x}^\star)$.
 Hence, C1 holds. Also, C2 holds, by assumptions of Theorem~\ref{theorem-asymptotic-normality-component-wise}. Now, we consider C3, which requires that the matrix $\boldsymbol{\Sigma} = a\,\mathbf{B} + \frac{1}{2}\mathbf{I}$ is stable, where $\mathbf{B} = -\phi^\prime(0) \nabla^2 f(\mathbf{x}^\star)$. Note that $\boldsymbol{\Sigma} = \frac{1}{2}\mathbf{I} - a\,\phi'(0) \nabla^2 f(\mathbf{x}^\star)$. Clearly, $\boldsymbol{\Sigma}$ is stable for a small enough $a$, because the matrix $\phi'(0) \nabla^2 f(\mathbf{x}^\star)$ is positive definite. More precisely, $\boldsymbol{\Sigma}$ is stable for $a>1/(2 \mu)$. Therefore, condition C3 holds, provided that $a>1/(2 \mu)$. We next consider condition C4. In the notation of Theorem~\ref{theorem-29}, consider the following quantity:
\begin{align}
    \mathbf{A}(t, \mathbf{x})  &\coloneqq \mathbb{E}\left[{\boldsymbol{\gamma}}(t+1,\mathbf{x},\omega) {\boldsymbol{\gamma}}(t+1,\mathbf{x},\omega)^\top \right] \\
    &= \mathbb{E} \left[ \left ( \phi(\nabla f(\mathbf{x})) - \Psi(\nabla f(\mathbf{x}) + \boldsymbol{\nu}^t)  \right )\left ( ( \phi(\nabla f(\mathbf{x})) - \Psi(\nabla f(\mathbf{x}) + \boldsymbol{\nu}^t) \right )^\top\right]
\end{align}
Now, because $\nabla f(\mathbf{x}^\star) = 0$; $\boldsymbol{\phi}(0) = 0$, and the entries of $\boldsymbol{\nu}^t$ are mutually independent, with pdf $p(u)$, we have that:
\begin{equation}
\label{eqn-mathcal-S}
\begin{split}
    \lim_{t\rightarrow \infty, \mathbf{x} \rightarrow \mathbf{x}^\star} \mathbf{A}(t,\mathbf{x}) =: \mathcal{S}_0 = \mathbb{E} \left [ \Psi(\boldsymbol{\nu}^t) \cdot \Psi(\boldsymbol{\nu}^t)^\top \right ] = \sigma_{\Psi}^2\cdot \mathbf{I},
    \end{split}
\end{equation}
where $\sigma_{\Psi}^2 = \int|\Psi(a)|^2p(a)da.$
\noindent 
Therefore, condition C4 holds. We finally verify condition C5. We follow the arguments analogous to those in Theorem~10 in \cite{kar2012distributed}. Condition~C5 means uniform integrability of the family $\{\|{\boldsymbol{\gamma}}(t+1,\mathbf{x},\omega)\|^2\}_{t =0,1,...,\,\|\mathbf{x}-\mathbf{x}^\star\|<\epsilon}$. 
\noindent 
We have:
\begin{equation}
    \|{\boldsymbol{\gamma}}(t+1,\mathbf{x},\omega)\|^2 \leq 2\|\phi(\nabla f(\mathbf{x}))\|^2 + 2\|\psi(\nabla f(\mathbf{x}) + \boldsymbol{\nu}^t)\|^2.
\end{equation}
We consider separately the cases when condition 2. or condition 3. hold in Assumption \ref{assumption-psi-comp-wise-iid-noise}. If condition 2. holds, then:
\begin{eqnarray}
    \|{\boldsymbol{\gamma}}(t+1,\mathbf{x},\omega)\|^2 &\leq& C_7 + C_8\|\mathbf{x}^t-\mathbf{x}^\star\|^2 + C_9\|\boldsymbol{\nu}^t\|^2\\
    &\leq & C_7 + C_8\,{\epsilon}^2 + C_9\|\boldsymbol{\nu}^t\|^2,
\end{eqnarray}
for some positive constants~$C_7,C_8,C_9$. Consider next the family
 $\{\widetilde{{\boldsymbol{\gamma}}}(t+1,\mathbf{x},\omega)\}_{t =0,1,...,\|\mathbf{x}-\mathbf{x}^\star\|<\epsilon}$,
 with
 \begin{equation}
  \widetilde{{\boldsymbol{\gamma}}}(t+1,\mathbf{x},\omega)   = C_7 + C_8\,{\epsilon}^2 + C_9\|\nu^t\|^2 .
 \end{equation}

The family  $\{\widetilde{{\boldsymbol{\gamma}}}(t+1,\mathbf{x},\omega)\}_{t=0,1,...,\|\mathbf{x}-\mathbf{x}^\star\|<\epsilon}$ is i.i.d. and hence it is uniformly integrable. The family $\{\|{\boldsymbol{\gamma}}(t+1,\mathbf{x},\omega)\|^2\}_{t =0,1,...,\|\mathbf{x}-\mathbf{x}^\star\|<\epsilon}$ is dominated by\\ $\{\widetilde{{\boldsymbol{\gamma}}}(t+1,\mathbf{x},\omega)\}_{t =0,1,...,\|\mathbf{x}-\mathbf{x}^\star\|<\epsilon}$ that is uniformly integrable, and hence \\$\{\|{\boldsymbol{\gamma}}(t+1,\mathbf{x},\omega)\|^2\}_{t =0,1,...,\|\mathbf{x}-\mathbf{x}^\star\|<\epsilon}$ is also uniformly integrable.
Hence, C5 holds under condition 2. of Assumption \ref{assumption-psi-comp-wise-iid-noise}. Now, let condition~3. in Assumption \ref{assumption-psi-comp-wise-iid-noise}) hold. Then:
\begin{eqnarray}
    \|{\boldsymbol{\gamma}}(t+1, \mathbf{x}, \omega)\|^2 &\leq& C_{10} + C_{11}\|\mathbf{x}-\mathbf{x}^\star\|^2\\
    &\leq& C_{10} + C_{11}\,{\epsilon}^2.
\end{eqnarray}
Consider the family
 $\{\widehat{{\boldsymbol{\gamma}}}(t+1,\mathbf{x},\omega)\}_{t =0,1,...,\|\mathbf{x}-\mathbf{x}^\star\|<\epsilon}$,
 with
 \begin{equation}
  \widehat{{\boldsymbol{\gamma}}}(t+1,\mathbf{x},\omega)   = C_{10} + C_{11}\,{\epsilon}^2.
 \end{equation}
The family  $\{\widehat{{\boldsymbol{\gamma}}}(t+1,\mathbf{x},\omega)\}_{t=0,1,...,\|\mathbf{x}-\mathbf{x}^\star\|<\epsilon}$ is uniformly integrable, and condition~C5 is verified analogously to the previous case.

 Summarizing, we have established that
 all conditions C1-C5
 of Theorem~\ref{theorem-29} hold true,
 thus the proof of  Theorem~\ref{theorem-asymptotic-normality-component-wise}. $\Box.$

\subsection{MSE analysis: Proof of Theorem~\ref{theorem-MSE-comp-Wise}}
%We next prove Theorem~\ref{theorem-MSE-comp-Wise}. 
We start with the following Lemma that shows that, with algorithm~\eqref{eqn-SGD-nonlinear}, almost surely, $\nabla f(\mathbf{x}^t)$ can be at most
$O(\mathrm{ln}(t))$.
\begin{lemma}
\label{Lemma-boundedness}
Let Assumptions \ref{assumption-f-i},
\ref{assumption-nu-comp-wise-iid-noise}, and Assumption~\ref{assumption-psi-comp-wise-iid-noise} with condition~3. hold (the nonlinearity with bounded outputs case). Then, for each $t =1,2,...$, we have:
 \begin{equation}
 \label{eqn-G-t}
 \| \nabla f(\mathbf{x}^t)\|
 \leq  G_t:=
 L \,\left( a\,C_2\,\sqrt{d}\,\frac{t^{1-\delta}}{1-\delta}
 +\|\mathbf{x}^0\|+\|\mathbf{x}^{\star}\|\right).
 \end{equation}
\end{lemma}

\textit{Proof}.
 Consider \eqref{eqn-SGD-nonlinear}.
 Because the output of each
 component nonlinearity $\Psi$ is
 bounded in the absolute value by $C_2$
  (Assumption~\ref{assumption-psi-comp-wise-iid-noise}), we have, for each $t \geq 1$:
  \begin{eqnarray}
  \nonumber
  \|\mathbf{x}^t\| &\leq&
   \|\mathbf{x}^0\| + a\,\sqrt{d}\,C_2\,\sum_{s=0}^{t-1} \frac{1}{(s+1)^{\delta}}\\
   &\leq&
   \label{eqn-lemma-proof-new}
   \|\mathbf{x}^0\|
   +
   a\,C_2\,\sqrt{d}\,\left( \frac{t^{1-\delta}}{1-\delta}\right).
  \end{eqnarray}
  Next, because $\nabla f$ is $L$-Lipschitz,
   we have:
   $
   \|\nabla f(\mathbf{x}^t)\| \leq
   L\,\|\mathbf{x}^t-\mathbf{x}^\star\|$. 
 Applying this inequality to \eqref{eqn-lemma-proof-new}, 
 the result follows. $\Box$

  We will also make use of the following Lemma.

 \begin{lemma}
 \label{lemma-G_t}
  There exists a positive constant $\xi$
   such that, for any $t=1,2,...$,
   there holds, almost surely, for each
   $j=1,...,d$, that:
   \[
   |\phi([\nabla f(\mathbf{x}^t)]_j)|
    \geq |[\nabla f(\mathbf{x}^t)]_j|\,
    \frac{\phi^\prime(0)\,\xi}{2\,G_t},
   \]
   where $G_t$
    is defined in~\eqref{eqn-G-t}.
 \end{lemma}

\textit{Proof}.
 Consider function~$\phi$ in~\eqref{eqn-phi-definition}.
 By Lemma~\ref{Lemma-Polyak},
 we have that $\phi^\prime(0)>0$
 and~$\phi$ is continuous at zero.\footnote{As $\phi$ is an odd function,
 for simplicity, in the proof we consider only
 nonnegative arguments of $\phi$, while analogous analysis applies for negative arguments of $\phi$.}
 Therefore, there exists a
 positive constant $\xi$ such that:
 \[
 \phi(a) \geq \frac{\phi^\prime(0)}{2}\,a,
 \]
 for any $a \in [0,\xi)$.
 Now, because $\phi$ is
 non-decreasing (by Lemma~\ref{Lemma-Polyak}),
  it holds for any $a^\prime>\xi$ that
  \begin{equation}
  \label{eqn-a-prime}
  \phi(a) \geq
  \frac{\phi^\prime(0)\,\xi\,a}{2\,a^\prime},
  \,\mathrm{for\,any\,}a\in [0,a^\prime).
  \end{equation} %
 Consider now $\nabla f(\mathbf{x}^t)$.
  By Lemma~\ref{Lemma-boundedness},
  we have that
  $\|\nabla f(\mathbf{x}^t)\| \leq G_t$, a.s.,
   and so, for any $j=1,...,d$,
   $|[\nabla f(\mathbf{x}^t)]_j|
   \leq G_t$.
   Therefore, in view of~\eqref{eqn-a-prime},
   setting $a^\prime=G_t$,
    the Lemma follows. $\Box$

We are now ready to prove Theorem~\ref{theorem-MSE-comp-Wise}.

\textit{Proof} (Proof of Theorem~\ref{theorem-MSE-comp-Wise}).
 Consider algorithm~\eqref{eqn-SGD-nonlinear}.
  By the Lipschitz property of $\nabla f$, we have, for any $\mathbf{x},
  \mathbf{y} \in {\mathbb R}^d$, that:
  \[
  f(\mathbf{y}) \leq f(\mathbf{x})
  +\nabla f(\mathbf{x})^\top 
  (\mathbf{y}-\mathbf{x})
  +\frac{L}{2}\|\mathbf{x}-\mathbf{y}\|^2,
  \]
  and so, 
  almost surely:
\begin{equation}
\label{eqn-1-proof}
\begin{split}
    f(\mathbf{x}^{t+1}) &\leq f(\mathbf{x}^t) + \bigl(\nabla f(\mathbf{x}^t)\bigr)^\top(-\alpha_t \Psi(\nabla f(\mathbf{x}^t) + \boldsymbol{\nu}^t))\\
    &+ \frac{L}{2} \alpha_t^2\|\Psi
    (\nabla f(\mathbf{x}^t) + \boldsymbol{\nu}^t)\|^2.
\end{split}
\end{equation}
 Next, letting
 $\boldsymbol{\eta}^t = \Psi
    (\nabla f(\mathbf{x}^t) + \boldsymbol{\nu}^t) -
    \boldsymbol{\phi}(\nabla f(\mathbf{x}^t))$,
    and using the fact that $\Psi$
     has bounded outputs, we obtain:
\begin{equation}
\label{eqn-2-proof}
\begin{split}
    f(\mathbf{x}^{t+1}) &\leq f(\mathbf{x}^t) + \bigl(\nabla f(\mathbf{x}^t)\bigr)^\top(-\alpha_t \boldsymbol{\phi}(\nabla f(\mathbf{x}^t) )) + \frac{L}{2} \alpha_t^2\,d^2 C_{12}^2 -\alpha_t
    \,(\nabla f(\mathbf{x}^t))^\top \boldsymbol{\eta}^t,\,\,\mathrm{a.s.},
\end{split}
\end{equation}
for some positive constant $C_{12}$. 
 Let $\mathcal{F}_t$
      be the history
      of~\eqref{eqn-SGD-nonlinear}
       up to iteration~$t$.
        Then, taking conditional expectation,
        and using the fact that
        $\mathbb{E}[\boldsymbol{\eta}^t\,|\,\mathcal{F}_t]=0$,
        we get that, almost surely:
\begin{equation}
\label{eqn-3-proof}
\begin{split}
    \mathbb{E}[f(\mathbf{x}^{t+1})\,|\,\mathcal{F}_t] &\leq f(\mathbf{x}^t) -\alpha_t\, \bigl(\nabla f(\mathbf{x}^t)\bigr)^\top\,\boldsymbol{\phi}(\nabla f(\mathbf{x}^t)) + \frac{L}{2} \alpha_t^2\,d^2\,C_{12}^2.
\end{split}
\end{equation}
  Next, using Lemma~\ref{lemma-G_t}, the fact that $\alpha_t=a/(t+1)^{\delta}$, and the fact that
  $G_t=O(t^{\epsilon})$, for $\epsilon>0$
  we obtain that. a.s.:
  \begin{equation}
\label{eqn-4-proof}
\begin{split}
    \mathbb{E}[f(\mathbf{x}^{t+1})\,|\,\mathcal{F}_t] &\leq f(x^t) -
    \frac{c^\prime}{(t+1)}\|\nabla f(\mathbf{x}^t)\|^2 + \frac{L}{2} \frac{a^2\,d^2\,C_2^2}{(t+1)^{2\delta}},
\end{split}
\end{equation}
        where $c^\prime =
        \frac{a\,(1-\delta)\xi\,
        \phi^{\prime}(0)}
  {2\,L \,(a\,C_2\,\sqrt{d}\,
 +\|\mathbf{x}^0\|+\|\mathbf{x}^{\star}\|)}$.
Next, by strong convexity of $f$,
we have that
$\|\nabla f(\mathbf{x}^t) - \nabla f(\mathbf{x}^\star)\|^2
 \geq 2\,\mu\,(f(\mathbf{x}^t) - f^\star)$.
Using the latter inequality, subtracting 
$f^\star$ from both sides of the inequality, taking expectation, and applying  Theorem~\ref{lemma-estimation}, claims (2) and (3), 
  we get the MSE rate result
  under condition 3. in Assumption~\ref{assumption-psi-joint-iid-noise}).

   We next consider the case when condition 2. in 
   Assumption~\ref{assumption-psi-joint-iid-noise}) holds, and we have the bounded second moment of $\boldsymbol{\nu}^t$.
    Following analogous arguments as in the first part of the proof,
    it can be shown that, a.s.:
\begin{equation}
\label{eqn-7-proof}
\begin{split}
    \mathbb{E}[f(\mathbf{x}^{t+1})\,|\,\mathcal{F}_t] &\leq f(\mathbf{x}^t) -
    \alpha_t\,\boldsymbol{\phi}(\nabla f(\mathbf{x}^t))^\top \nabla f(\mathbf{x}^t) + \frac{L}{2} \alpha_t^2
    \left( C_{13} + C_{14} \|\boldsymbol{\nu}^t\|^2\right),
\end{split}
\end{equation}
for some positive constants $C_{13},C_{14}$. 
Next, because $\inf_{a\neq 0}\frac{|\phi(a)|}{|a|}>0$,
we have that
$\phi(\nabla f(\mathbf{x}^t))^\top \nabla f(\mathbf{x}^t)
\geq C_{15}\,\|\nabla f(\mathbf{x}^t) \|^2$, for some constant  $C_{15}>0$.
 Using the latter bound in \eqref{eqn-7-proof}, subtracting 
 $f^\star$ from both sides of the inequality, 
  taking expectation, and
 applying Theorem~\ref{lemma-estimation}, 
 claim (1) and (2), 
  the result follows.
   $\Box$

\section{Intermediate results and proofs: Joint nonlinearities}
\label{section-proofs-joint}
Subsection~{6.1} provides the required intermediate results, 
while Subsection~{6.2} proves Theorem~\ref{theorem-MSE-joint-iid}.

\subsection{Intermediate results: Joint nonlinearities}
Recall function $\mathcal{N}:\mathbb{R}_+ \mapsto \mathbb{R}_+$ 
in Assumption~\ref{assumption-psi-joint-iid-noise}. 
We first state and prove the following Lemma
on the properties of function~$\mathcal{N}$.
 \begin{lemma}
 \label{lemma2.1}
 Under Assumption \ref{assumption-psi-joint-iid-noise},
 for any $\mathbf{x},
 \mathbf{u}\in {\mathbb R}^d$, 
 such that $\|\mathbf{u}\| > \|\mathbf{x}\| $, there holds:
 \begin{equation}
 \begin{split}
     |\mathcal{N}(\|
     \mathbf{x}+\mathbf{u}\|) 
     &- \mathcal{N}(\|\mathbf{x}-\mathbf{u}\|)| \leq \frac{\|\mathbf{x}\|}{\|\mathbf{u}\|}
     \left[ \mathcal{N}(\|\mathbf{x}+\mathbf{u}\|) + \mathcal{N}(\|\mathbf{x}-\mathbf{u}\|) \right].
     \label{lemma31}
\end{split}
 \end{equation}
 \end{lemma}
 \textit{Proof.} Fix a pair $\mathbf{x},\mathbf{u}\in {\mathbb R}^d$, such that 
 $\|\mathbf{u}\| > \|\mathbf{x}\|$, and assume without loss of generality that $\mathcal{N}(\|\mathbf{x}+\mathbf{u}\|) \geq \mathcal{N}(\|\mathbf{x}-\mathbf{u}\|)$.
Then, \eqref{lemma31}  is equivalent to:
  \begin{equation}
     (\|\mathbf{u}\| - \|\mathbf{x}\|) \mathcal{N}(\|\mathbf{x}+\mathbf{u}\|)
     \leq ( \|\mathbf{u}\| + \|\mathbf{x}\|) \mathcal{N}(\|\mathbf{x}-\mathbf{u}\|).
 \label{lemma31eq2}
 \end{equation}
  Denote by $\rho=\|\mathbf{u}\|$. Notice that:
 \begin{equation}
 \label{eqn-proof-aux}
     \rho - \|\mathbf{x}\| 
     \leq \|\mathbf{x} + \mathbf{u}\| 
     \leq \|\mathbf{x}\| + \|\mathbf{u}\| 
     = \|\mathbf{x}\| + \rho,
 \end{equation}
 and similarly,
 \begin{equation*}
     \rho + \|\mathbf{x}\| 
     \geq \|\mathbf{x} - \mathbf{u}\| \geq \rho - \|\mathbf{x}\|.
 \end{equation*}
 As $\mathcal{N}$ is non-increasing, it follows that:
 \begin{equation*}
     \mathcal{N}(\|\mathbf{x}+\mathbf{u}\|) \leq \mathcal{N}(\rho - \|\mathbf{x}\|),\,\,\,
     \mathcal{N}(\|\mathbf{x}-\mathbf{u}\|) \geq \mathcal{N}(\rho + \|\mathbf{x}\|).
 \end{equation*}
 Now, we have:
 \begin{equation}
 \label{eqn-proof-C1}
     (\|\mathbf{u}\| - \|\mathbf{x}\|) \mathcal{N}(\|\mathbf{x}+\mathbf{u}\|) \leq (\rho - \|\mathbf{x}\|) 
     \mathcal{N}(\rho - \|\mathbf{x}\|),
 \end{equation}
 and similarly:
 \begin{equation}
  \label{eqn-proof-C2}
     ( \|\mathbf{u}\| + \|\mathbf{x}\|) \mathcal{N}(\|\mathbf{x}-\mathbf{u}\|) \geq (\rho + \|\mathbf{x}\|) 
     \mathcal{N}(\rho + \|\mathbf{x}\|).
 \end{equation}
 By assumption, function $a\mapsto a \mathcal{N}(a), \, a>0$,  
 is non-decreasing, and so 
 $( \rho - \|\mathbf{x}\|) \mathcal{N}(\rho-\|\mathbf{x}\|)$
  $\leq ( \rho + \|\mathbf{x}\|) 
  \mathcal{N}(\|\mathbf{x}\|+\rho)$. Thus, combining~\eqref{eqn-proof-C1}
 and~\eqref{eqn-proof-C2}, we have that \eqref{lemma31eq2} holds, which is in turn equivalent to the claim of the Lemma.

We now define map~${\boldsymbol{\phi}} :\, \mathbb{R}^d \mapsto \mathbb{R}^d$, as follows. For a fixed (deterministic) point $\mathbf{w} \in {\mathbb R}^d$, we let:
 \begin{equation}
 \label{eqn-phi-joint}
     {\boldsymbol{\phi}}(\mathbf{w}) = \int \boldsymbol{\Psi}(\mathbf{w}+\mathbf{u})p(\mathbf{u})d\mathbf{u} = \mathbb{E}[\boldsymbol{\Psi}(\mathbf{w}+\boldsymbol{\nu}^0)],
 \end{equation}
 where the expectation is taken with respect to the joint pdf
 of the gradient noise at any iteration~$t$, e.g., $t=0$.
 The map ${\boldsymbol{\phi}} :\, \mathbb{R}^d \mapsto \mathbb{R}^d$ is, abusing notation, a counterpart of the
 component-wise map~$\phi :\, \mathbb{R} \mapsto \mathbb{R}$ in~\eqref{eqn-phi-definition}. We have the following Lemma.
 
\begin{lemma}
\label{lemma2.2}
 The following holds: 
 \begin{equation}
     {\boldsymbol{\phi}}(\mathbf{x})^\top \mathbf{x} \geq 
     2(1-{\kappa})
     \|\mathbf{x}\|^2
     \int\displaylimits_{{\mathcal{J}}(\mathbf{x})}
     \mathcal{N}(\|\mathbf{x}\|+\|\mathbf{u}\|)
     p(\mathbf{u})d\mathbf{u},
 \end{equation}
 where ${\mathcal{J}}(\mathbf{x}) = \{\mathbf{u}: \,
 \frac{\mathbf{u}^\top \mathbf{x}}{\|\mathbf{u}\|\|\mathbf{x}\|} \in [0,{\kappa}]\}$, and ${\kappa}$  is any constant 
 in the interval~$(0,1)$.
\end{lemma}
%
%The Lemma provides a 
\textit{Proof.} Let us fix arbitrary $\mathbf{x}\in \mathbb{R}^d, \mathbf{x}\neq0$. As $\boldsymbol{\Psi}(\mathbf{a}) = \mathbf{a}\mathcal{N}(\|\mathbf{a}\|)$, we have:
 \begin{align}
     {\boldsymbol{\phi}}(\mathbf{x})^\top \mathbf{x} 
     &= \int_{\mathbf{u}\in \mathbb{R}^d}
     \underbrace{(\mathbf{x}+\mathbf{u})^\top \mathbf{x} \mathcal{N}(\|\mathbf{x}+\mathbf{u}\|)}_{\coloneqq \mathcal{M}(\mathbf{x},\mathbf{u})}
     p(\mathbf{u})d\mathbf{u} \\ &= \int_{J_1(\mathbf{x})=\{\mathbf{u}:\,
     \mathbf{u}^\top \mathbf{x} \geq 0 \}}\mathcal{M}(\mathbf{x},\mathbf{u})
     p(\mathbf{u})d\mathbf{u} \\ 
     &+ \int_{J_2(\mathbf{x})=
     \{\mathbf{u}:\,\mathbf{u}^\top \mathbf{x} < 0 \}}\mathcal{M}(\mathbf{x},\mathbf{u})
     p(\mathbf{u})d\mathbf{u}.
 \end{align}
Note also that there holds:
 \begin{equation*}
     \mathcal{M}(\mathbf{x},
     \mathbf{u}) = 
     (\|\mathbf{x}\|^2 + \mathbf{u}^\top \mathbf{x})\mathcal{N}(\|\mathbf{x}+\mathbf{u}\|).
 \end{equation*}
 Similarly, 
 \begin{equation*}
     \mathcal{M}(\mathbf{x},-\mathbf{u}) 
     = (\|\mathbf{x}\|^2 - \mathbf{u}^\top \mathbf{x})\mathcal{N}(\|\mathbf{x}-
     \mathbf{u}\|).
 \end{equation*}
 Therefore, using the fact that 
$p(\mathbf{u}) = p(-\mathbf{u}),$ 
for all $u \in \mathbb{R}^d$, we obtain:
 \begin{equation}
 {\boldsymbol{\phi}}(\mathbf{x})^\top \mathbf{x} = \int_{J_1(\mathbf{x})}\mathcal{M}_2(\mathbf{x},
 \mathbf{u})p(\mathbf{u})d\mathbf{u},
 \end{equation}
 where 
 $\mathcal{M}_2(\mathbf{x},\mathbf{u}) 
 = [(\|\mathbf{x}\|^2 + \mathbf{u}^\top \mathbf{x})
 \mathcal{N}(\|\mathbf{x}+\mathbf{u}\|) 
 + (\|\mathbf{x}\|^2 - \mathbf{u}^\top \mathbf{x})
 \mathcal{N}(\|\mathbf{x}-\mathbf{u}\|)]$. There holds:
 \begin{equation}
 \begin{split}
     \mathcal{M}_2(\mathbf{x},\mathbf{u}) &\geq \|\mathbf{x}\|^2
     [\mathcal{N}(\|\mathbf{x}+\mathbf{u}\|) + \mathcal{N}(\|\mathbf{x}-\mathbf{u}\|)] - \|\mathbf{u}\|\|\mathbf{x}\|
     |\mathcal{N}(\|\mathbf{x}+\mathbf{u}\|) - \mathcal{N}(\|\mathbf{x}-\mathbf{u}\|)|.
 \end{split}
 \end{equation}
 Since $\mathbf{u}\in J_1(\mathbf{x})$, there holds $\|\mathbf{x}+\mathbf{u}\|
 \geq\|\mathbf{x}-\mathbf{u}\|$. Now, using Lemma \ref{lemma2.1}, we have:
 \begin{equation}
 \begin{split}
     \mathcal{M}_2(\mathbf{x},\mathbf{u}) &\geq \|\mathbf{x}\|^2
     [\mathcal{N}(\|\mathbf{x}+\mathbf{u}\|) + \mathcal{N}(\|\mathbf{x}-\mathbf{u}\|)] 
     - \phantom{-} \|\mathbf{u}\|\|\mathbf{x}\|
     \frac{\|\mathbf{x}\|}{\|\mathbf{u}\|}|
     \mathcal{N}(\|\mathbf{x}+\mathbf{u}\|) + \mathcal{N}(\|\mathbf{x}-\mathbf{u}\|)| = 0.
 \end{split}
 \end{equation}
 Therefore, we have:
 \begin{equation}
     \mathcal{M}_2(\mathbf{x},\mathbf{u}) 
     \geq 0, \, \mathrm{for\,any}\, 
     \mathbf{u} \in J_1(\mathbf{x}), \,\|\mathbf{u}\| > \|\mathbf{x}\|.
 \end{equation}
 Now, 
 consider ${\mathcal{J}}(\mathbf{x}) = 
 \{\mathbf{u}\in \mathbb{R}^d: \,
 \mathbf{u}^\top \mathbf{x} \geq 0, \,
 \frac{\mathbf{u}^\top \mathbf{x}}
 {\|\mathbf{u}\|\|\mathbf{x}\|} \in [0, {\kappa}]\}$, where ${\kappa} \in (0, 1)$. 
 Let us consider $\mathbf{u} \in {\mathcal{J}}(\mathbf{x})$ such that $\|\mathbf{u}\| > 
 \|\mathbf{x}\|$. Then, using Lemma~{6.1}, we get:
\begin{equation}
 \begin{split}
     \mathcal{M}_2(\mathbf{x},\mathbf{u}) &\geq \|\mathbf{x}\|^2
     [\mathcal{N}(\|\mathbf{x}+\mathbf{u}\|) + \mathcal{N}(\|\mathbf{x}-\mathbf{u}\|)] 
     \\
     &- \|\mathbf{u}\|
     \|\mathbf{x}\|{\kappa}\underbrace{|\mathcal{N}(\|\mathbf{x}+\mathbf{u}\|) - \mathcal{N}(\|\mathbf{x}-\mathbf{u}\|)|}
      \\
     &\geq (1 - {\kappa})
     \|\mathbf{x}\|^2 ( \mathcal{N}(\|\mathbf{x}+\mathbf{u}\|) + \mathcal{N}(\|\mathbf{x}-\mathbf{u}\|)).
 \end{split}
 \end{equation}
 Now, consider $\mathbf{u} \in {\mathcal{J}}(\mathbf{x})$ such that $\|\mathbf{u}\| \leq \|\mathbf{x}\|$. Then, there holds:
 \begin{equation}
 \begin{split}
     \mathcal{M}_2(\mathbf{x},\mathbf{u}) &\geq \|\mathbf{x}\|^2[\mathcal{N}(\|\mathbf{x}+
     \mathbf{u}\|) 
     + \mathcal{N}(\|\mathbf{x}-\mathbf{u}\|)] - \\
     &\phantom{-} 
     \underbrace{\|\mathbf{u}\|}_{\leq \|\mathbf{x}\|}\|\mathbf{x}\|{\kappa}
     |\underbrace{\mathcal{N}(\|\mathbf{x}+
     \mathbf{u}\|)}_{\geq 0} + \underbrace{\mathcal{N}(\|\mathbf{x}-\mathbf{u}\|)|}_{\geq 0} \\
     &\geq (1-{\kappa})\|\mathbf{x}\|^2
     (\mathcal{N}(\|\mathbf{x}+\mathbf{u}\|) + \mathcal{N}(\|\mathbf{x}-\mathbf{u}\|)).
 \end{split}
 \end{equation}
 where the last inequality holds due to the fact that $|a-b| \leq |a|+|b|$, for any $a,b \in \mathbb R$.  Now, we have:
 \begin{equation}
 \begin{split}
     \mathcal{M}_2(\mathbf{x},\mathbf{u}) &\geq (1-{\kappa})
     \|\mathbf{x}\|^2
     (\underbrace{\mathcal{N}
     (\|\mathbf{x}+\mathbf{u}\|)}_{\geq\mathcal{N}(\|x\| + \|u\|)} + \underbrace{\mathcal{N}
     (\|\mathbf{x}-\mathbf{u}\|)}_{\geq\mathcal{N}(\|\mathbf{x}\| + \|\mathbf{u}\|)}) \\
     &\geq 2(1-{\kappa})\|\mathbf{x}\|^2\mathcal{N}(\|\mathbf{x}\|+
     \|\mathbf{u}\|), \,\mathrm{for\,any}\, \mathbf{u} \in 
     {\mathcal{J}}(\mathbf{x}).
 \end{split}
 \end{equation}
 Combining (6.15) and (6.17), 
 we finally get:
 \begin{equation}
 \begin{split}
     {\boldsymbol{\phi}}(\mathbf{x})^\top \mathbf{x} &\geq \int\displaylimits_{{\mathcal{J}}(\mathbf{x})}
     2(1-{\kappa})
     \|\mathbf{x}\|^2\mathcal{N}(\|\mathbf{x}\| 
     +\|\mathbf{u}\|)p(\mathbf{u})d\mathbf{u} \\
     &=2(1-{\kappa})
     \|\mathbf{x}\|^2
     \int\displaylimits_{{\mathcal{J}}(\mathbf{x})}
     \mathcal{N}(\|\mathbf{x}\|+\|\mathbf{u}\|)
     p(\mathbf{u})d\mathbf{u}.
 \end{split}
 \end{equation}

\begin{lemma}
\label{Lemma-boundedness-joint}
Let Assumptions \ref{assumption-f-i},
\ref{assumption-p-joint-iid-noise}, and  Assumption~\ref{assumption-psi-joint-iid-noise} with condition~3. hold (the nonlinearity with bounded outputs case). Then, for each $t =1,2,...$, we have:
 \begin{equation}
 \label{eqn-G-t-joint}
 \| \nabla f(\mathbf{x}^t)\|
 \leq  
 G_t^{\prime}:=
 L \,(\,a\,C_2^{\prime}\,\frac{t^{1-\delta}}{1-\delta}
 +\|\mathbf{x}^0\|+\|\mathbf{x}^{\star}\|\,).
 \end{equation}
\end{lemma}

\textit{Proof.} The proof is analogous to the proof of Lemma~\ref{Lemma-boundedness}.

\subsection{Proof of Theorem~\ref{theorem-MSE-joint-iid}: Joint nonlinearities}
\textit{Proof.} We first consider the case with bounded nonlinearity (condition 3. in Assumption~\ref{assumption-psi-joint-iid-noise}). Analogously to the proof of \ref{theorem-MSE-comp-Wise}, it can be shown that, a.s.:
\begin{equation}
\label{eqn-8-proof}
\begin{split}
    \mathbb{E}[f(\mathbf{x}^{t+1})\,|\,
    \mathcal{F}_t] &\leq f(\mathbf{x}^t) -
    \alpha_t\,{\boldsymbol{\phi}}(\nabla f(\mathbf{x}^t))^\top \nabla f(\mathbf{x}^t) + \alpha_t^2 \,C_{17},
\end{split}
\end{equation}
for some positive constant~$C_{17}$. By Lemma \ref{lemma2.2}, there holds, for
$\mathbf{a}:=\nabla f(\mathbf{x}^t)$, a.s.:
\begin{equation}
\label{eqn-lower-bound-N-joint}
    \bigl( {\boldsymbol{\phi}}(\mathbf{a}) \bigr)^\top \mathbf{a} \geq 2 (1-\kappa)\|\mathbf{a}\|^2\int_{\mathcal{J}}\mathcal{N}
    (\|\mathbf{a}\| + \|\mathbf{u}\|)
    p(\mathbf{u})d\mathbf{u},
\end{equation}
where we recall $\mathcal{J} = \{\mathbf{u}:\,\frac{\mathbf{u}^\top \mathbf{a}}{\|\mathbf{u}\|\|\mathbf{a}\|}
\in[0,\kappa]\}$, where $\kappa \in (0,1)$ is a constant.  
Note that, as $a \mapsto a\,\mathcal{N}(a)$ is non-decreasing, $\mathcal{N}$ satisfies: $\mathcal{N}(b) \geq \mathrm{min} \left(\frac{\mathcal{N}(1)}{b}, {\mathcal{N}(1)}\right)\,\mathrm{for\,any}\, b>0$. Consider constant $B_0$ in condition 2. of Assumption~4. Then, for all $\mathbf{u}$ such that $\|\mathbf{u}\|\leq B_0$, there holds $\mathcal{N}(\|\mathbf{a}\| + \|\mathbf{u}\|) \geq \mathrm{min} \left(\frac{\mathcal{N}(1)}{\|\mathbf{a}\| + B_0}, {\mathcal{N}(1)}\right)$. Therefore, we have that, almost surely, for sufficiently large $t$:
\begin{eqnarray*}
&\,& \|\nabla f(\mathbf{x}^t)\|^2\int_{J_4}
\mathcal{N}(\|\nabla f(\mathbf{x}^t)\|+\|\mathbf{u}\|)p(\mathbf{u})d
\mathbf{u} \geq C_{18}\,\frac{\|\nabla f(\mathbf{x}^t)\|^2}{
G_t^{\prime} + B_0},
\end{eqnarray*}
for some positive constant~$C_{18}$. Here, $J_4 = \{u\in \mathbb{R}^d: 
\frac{\mathbf{u}^\top \nabla f(\mathbf{x}^t)}{\|\mathbf{u}\|\|\nabla f(\mathbf{x}^t)\|} \in [0, {\kappa}],\,\|\mathbf{u}\| \leq B_0\}$. Combining the last bound with Lemmas~\ref{lemma2.2} and~\ref{Lemma-boundedness-joint}, in view of condition 2. in Assumption~4, we obtain that, for sufficiently large $t$, a.s.:
 
\begin{equation}
    \label{eqn-proof-MSE-rate-2}
({\boldsymbol{\phi}}(\nabla f(\mathbf{x}^t)))^\top
\nabla f(\mathbf{x}^t)
\geq C_{19}\,\frac{\|\nabla f(\mathbf{x}^t)\|^2}{B_0 + G_t^{\prime}},
\end{equation}
where the positive constant~$C_{19}$
 can be taken as
 $C_{19} =2 (1-\kappa)\lambda(\kappa)\mathcal{N}(1)$. 
 
 Applying the bound \eqref{eqn-proof-MSE-rate-2} to \eqref{eqn-8-proof} we obtain an equivalent
 to~\ref{eqn-7-proof}. Therein, $c^\prime$ in~\ref{eqn-7-proof} is replaced with a positive constant $c^{\prime \prime}$ that can be taken as $c^{\prime \prime} = \frac{2\,a\,(1-\kappa)\lambda(\kappa)(1-\delta)\mathcal{N}(1)}{L\,(\,a\,C_2^\prime +\|\mathbf{x}^0\| +\|\mathbf{x}^\star\|\,)+B_0)}$. We now proceed analogously to the proof of Theorem~\ref{theorem-MSE-comp-Wise}, by applying claims (2) and (3) of Theorem~\ref{lemma-estimation}. The result
for the bounded nonlinearity $\boldsymbol{\Psi}$ follows, with the rate $\zeta$ being any positive number less than 
\begin{equation}
\label{eqn-zeta-joint}
\mathrm{min} \left\{ 
2 \delta-1, \, \frac{2\,a\,\mu\,(1-\kappa)\lambda(\kappa)(1-\delta)\mathcal{N}(1)}{L\,(\,a\,C_2^\prime +\|\mathbf{x}^0\|
  +\|\mathbf{x}^\star\|\,)+B_0}
\right\}.
\end{equation}

We now prove the alternative case, for the nonlinearity $\boldsymbol{\Psi}$ with unbounded outputs and finite second moment of $\boldsymbol{\nu}^t$. We have that  $\inf_{\mathbf{x} \neq 0}\frac{\|\boldsymbol{\Psi}(\mathbf{x})\|}{\|\mathbf{x}\|}>0$. This is equivalent to saying that $\mathcal{N}$ is lower-bounded by a positive constant, i.e., $\mathcal{N}(a) \geq C_{20}$, for each $a$, for some constant $C_{20}>0$. Then, it follows that, a.s.:
\begin{equation}
({\boldsymbol{\phi}}(\nabla f(\mathbf{x}^t)))^\top
\nabla f(\mathbf{x}^t)
\geq C_{21}\,
 \|\nabla f(\mathbf{x}^t)\|^2,
\end{equation}
for some positive constant~$C_{21}$.  The proof then proceeds analogously to the proof of Theorem~\ref{theorem-MSE-comp-Wise} by applying the appropriate variant of Theorem~\ref{lemma-estimation}.

It remains to prove a.s. convergence of \eqref{eqn-SGD-nonlinear}. We do so again by verifying conditions B1-B5 in Theorem~{8.1}. Algorithm \eqref{eqn-SGD-nonlinear} admits again the representation in Theorem~{8.1} with 
 \begin{eqnarray}
\label{eqn-R-new}
    &\,& {\mathbf{r}}(\mathbf{x}) 
    = -\boldsymbol{{\boldsymbol{\phi}}}(\nabla f(\mathbf{x}))\\
    &\,& 
    {\boldsymbol{\gamma}}(t+1, \mathbf{x}, \omega) = 
    \boldsymbol{{\boldsymbol{\phi}}}(\nabla f(\mathbf{x})) - \boldsymbol{\Psi}(\nabla f(\mathbf{x}) + \boldsymbol{\nu}^t).
\end{eqnarray}
Conditions B1, B2 clearly hold. Condition B3 follows from Lemma~\ref{lemma2.2}. Condition B4 holds analogously to the proof of Theorem~\ref{theorem-almost-sure-component-wise}. Finally, condition B5 follows from the definition of the step-size sequence in Theorem~\ref{theorem-MSE-joint-iid}. Thus, the result.
$\Box$

\section{Experiments}
\label{section-experiments}
In order to benchmark the proposed nonlinear SGD framework, we consider \texttt{Heart}, \texttt{Diabetes} and \texttt{Australian} datasets from the LibSVM library~\cite{chang2011libsvm}. We consider the logistic regression loss function for binary classification, see, e.g.,~\cite{gorbunov2020stochastic}, where function~$f$ in Eq.~\eqref{eq:opt_problem} is the empirical loss, i.e., the sum of the logistic losses across all data points in a given dataset. 

As it has been studied in~\cite{gorbunov2020stochastic} (see Figure~2 in~\cite{gorbunov2020stochastic}), we have, near the solution~$\mathbf{x}^\star$, the following behavior with respect to gradient noise. (See also~\cite{gorbunov2020stochastic} for details how the gradient noise is evaluated in Figure~2 therein.)  With the \textsc{heart} dataset, tails of stochastic gradients are not heavy. On the other hand, for \textsc{diabetes} and \textsc{australian} datasets, the gradient noise has outliers and exhibits a heavy-tail behavior. 

We consider three different nonlinearities to demonstrate the effectiveness of our nonlinear framework, namely, \texttt{tanh} (hyperbolic tangent), \texttt{sign} and a bi-level customization of \texttt{sign} with ${\Psi}(x)=-1,$ $ -0.5, 0.5, 1,$ for $x\in(-\infty,-0.5], (-0.5, 0],$ $ (0, 0.5], (0.5, \infty]$,~respectively (\texttt{nonlinear-quantizer} in figures). Note that the  \texttt{tanh} function may be considered a smooth approximation of \texttt{sign}. We benchmark the above methods against the linear SGD, clipped-SGD and SSTM along with a clipped version of SSTM from \cite{gorbunov2020stochastic}. For each of the methods, we use batch sizes of $50$, $100$ and $20$  for the \texttt{Australian}, \texttt{Diabetes} and \texttt{Heart} datasets,  respectively. We also consider clipped-SGD with periodically decreasing clipping level (\texttt{d-clipped-SGD} in Figures) as a baseline as introduced in \cite{gorbunov2020stochastic}. This method starts with some initial clipping level and after every $l$ epochs the clipping level is multiplied by some constant $c\in (0,1)$. The step sizes $\alpha_t$ (learning rates) for each method from our framework were tuned after an experimentation. The learning rates for the baselines, i.e., SGD, clipped-SGD, SSTM and clipped-SSTM are also tuned and are selected to be as in \cite{gorbunov2020stochastic}. In more detail, the learning rates for the proposed methods are of the form $a/(b\,(t+1)+L)$, where we recall that $t$ is the iteration counter, $L$ is the smoothness constant of~$\nabla f$, and parameters~$a,b$ are tuned via grid search. The value of $a$ is chosen to be $1.0$, $1.5$ and $5.0$, respectively, for \texttt{Heart}, \texttt{Diabetes} and \texttt{Australian} and for all the three non-linearities. The value of $b$ is chosen to be $0.001$, $7.0$ and $7.0$ respectively for \texttt{Australian}, \texttt{Heart} and \texttt{Diabetes} datasets for the \texttt{sign} nonlinearity. The value of $b$ is chosen to be $0.0001$, $2.0$ and $3.0\times 10^{-6}$ respectively for \texttt{Australian}, \texttt{Heart} and \texttt{Diabetes} datasets for the \texttt{tanh} nonlinearity. The value of $b$ is chosen to be $0.001$, $5.0$ and $5.0$ respectively for \texttt{Australian}, \texttt{Heart} and \texttt{Diabetes} datasets for the \texttt{nonlinear-quantizer} nonlinearity. 

We first note that (see Figure \ref{fig:heavy_tail_comp}) \texttt{d-clipped-SGD} stabilizes the trajectory as compared to the linear SGD, even if the initial clipping level was high. At the same time, clipped-SGD with large clipping levels performs similarly as SGD. It is noteworthy, that SGD has the least oscillations for \texttt{Australian} and \texttt{Diabetes} datasets, despite the fact that these datasets have heavier or similar tails. This can be attributed to the fact that SGD does not get close to the solution in terms of functional value.  SSTM in particular shows large oscillations, which can be attributed to it being a version of accelerated/momentum-based methods and its usage of small batch sizes. \texttt{Clipped-SSTM} on the other hand suffers less from  oscillations and has a comparable convergence rate as SSTM. In comparison, all the three nonlinear schemes that have been proposed in this paper, have very little oscillations. While the \texttt{tanh} algorithm is outperformed by the algorithms with other nonlinearities from our framework, its performance is at par with the other baselines from \cite{gorbunov2020stochastic}. In particular, the \texttt{sign} algorithm compares favorably to  other baselines in terms of convergence for \texttt{Australian} and \texttt{Heart} datasets. The  \texttt{nonlinear-quantizer} algorithm outperforms other baselines for the \texttt{Diabetes} dataset. The good behavior of \texttt{tanh} and \texttt{sign} on the heavy-tail data sets, specially relative to the linear SGD, also viewing \texttt{tanh} as a smooth approximation of \texttt{sign}, might also be related with the insights from Example~{3.4}. In summary, the three simple example nonlinearities from the proposed framework are comparable or favorable over the considered state of the art benchmarks on the studied datasets. 
\vspace{-2mm}

\section{Conclusion}
\label{section-conclusion}
We proposed a general framework for nonlinear stochastic gradient descent (SGD) under heavy-tail gradient noise. Unlike existing studies of SGD under heavy-tail noise that focus on specific nonlinear functions (e.g., adaptive clipping), our framework includes a broad class of component-wise (e.g., sign gradient) and joint (e.g., gradient clipping) nonlinearities. We establish for the considered methods almost sure convergence, MSE convergence rate, and also asymptotic covariance for component-wise nonlinearities. We carry out numerical experiments on several real datasets that exhibit heavy tail gradient noise effects. The experiments show that, while our framework is more general than existing studies of SGD under heavy-tail noise, several  easy-to-implement nonlinearities from our framework are competitive with state of the art alternatives.

\section*{Appendix}
\label{section-appendix}

\subsection*{\textbf{A}. Some results in stochastic approximation}
\label{section-appendixA}

We present a useful result on single time scale stochastic
approximation; see~\cite{nevelson1976stochastic}, Theorems 4.4.4 and 6.6.1.  
\begin{theorem}
\label{theorem-29}
Let $\left\{\mathbf{x}^t\in\mathbb{R}^{d}\right\}$ be a random  sequence
that satisfies:
\begin{equation}
\label{RM:1}
\mathbf{x}^{t+1}=\mathbf{x}^t+\alpha_t\left[{\mathbf{r}}(\mathbf{x}^t)+{\boldsymbol{\gamma}}
\left(t+1,\mathbf{x}^t,\omega\right)\right],
\end{equation}
where, ${\mathbf{r}}(\cdot):\mathbb{R}^{d}\longmapsto\mathbb{R}^{d}$ is Borel measurable and
$\left\{{\boldsymbol{\gamma}}(t,\mathbf{x},\omega)\right\}_{t\geq
0,~\mathbf{x}\in\mathbb{R}^{d}}$ is a family of random
vectors in $\mathbb{R}^{d}$, defined on a probability
space $(\Omega,\mathcal{F},\mathcal{P})$, and $\omega\in\Omega$ is
a canonical element. Let the following sets of
assumptions hold:
\begin{itemize}
\item{\textbf{(B1)}}: The function ${\boldsymbol{\gamma}}(t,\cdot,
\cdot):\mathbb{R}^{d}\times\Omega\longrightarrow\mathbb{R}^{d}$ is
$\mathcal{B}^{d}\otimes\mathcal{F}$
measurable for every $t$; $\mathcal{B}^{d}$ is the Borel algebra of $\mathbb{R}^{d}$.
\item{\textbf{(B2)}}: There exists a filtration
$\left\{\mathcal{F}_{t}\right\}_{t\geq 0}$ of $\mathcal{F}$, such that, for
each $t$, the family of random vectors
$\left\{{\boldsymbol{\gamma}}\left(t,\mathbf{x},\omega\right)
\right\}_{\mathbf{x}\in\mathbb{R}^{d}}$ is $\mathcal{F}_{t}$ measurable, zero-mean and independent of
$\mathcal{F}_{t-1}$.
\item{\textbf{(B3)}}: There exists a twice continuously differentiable function
$V \left( \mathbf{x} \right) $ with bounded second order partial derivatives and a point $\mathbf{x}^{{\star}}\in\mathbb{R}^{d}$ satisfying:
\begin{align*}
&V\left(\mathbf{x}^{{\star}}\right)=0,
\:\:V\left(\mathbf{x}\right)>0,\,\mathbf{x}\neq\mathbf{x}^{{\star}},
\:\:\lim_{\left\|\mathbf{x}\right\|\rightarrow\infty}V\left(\mathbf{x}\right)=\infty,
\\
&
\sup_{\epsilon<\left\|\mathbf{x}-\mathbf{x}^{{\star}}\right\|<
\frac{1}{\epsilon}}\left(\mathbf{r}\left(\mathbf{x}\right),
V_{\mathbf{x}}\left(\mathbf{x}\right)\right)<0,\:\mathrm{for\,any}\,\epsilon>0
%&
\end{align*}
where $V_{\mathbf{x}}\left(\mathbf{x}\right)$ denotes the gradient (vector) of $V(\cdot)$ at $\mathbf{x}$.

\item{\textbf{(B4)}}: There exist constants $k_{1},k_{2}>0$, such
that,
\begin{align*}
&
\hspace{-.6cm}
\left\|\mathbf{r}\left(\mathbf{x}\right)\right\|^{2}+\mathbb{E}\left[\left\|
{\boldsymbol{\gamma}}\left(t+1,\mathbf{x},\omega\right)\right\|^{2}\right]
\leq
k_{1}\left(1+V\left(\mathbf{x}\right)\right)-
\\
&
\hspace{3.5cm}
-
k_{2}\left(\mathbf{r}\left(\mathbf{x}\right),
V_{\mathbf{x}}\left(\mathbf{x}\right)\right)
\end{align*}
\item{\textbf{(B5)}}: The weight sequence $\left\{\alpha_t\right\}$ satisfies
\begin{equation}
\label{alphacond-b} \alpha_t>0,\:\sum_{t\geq 0}\alpha_t=\infty, \: \sum_{t\geq 0}\alpha_t^{2}<\infty.
\end{equation}
\end{itemize}

\begin{itemize}

\item{\textbf{(C1)}}: The function $\mathbf{r}\left(\mathbf{x}\right)$ admits the
representation
\begin{equation}
\label{assc1.1}
\mathbf{r}\left(\mathbf{x}\right)=B\left(\mathbf{x}-\mathbf{x}^{{\star}}\right)
+\delta\left(\mathbf{x}\right)
\end{equation}
where
\begin{equation}
\label{assc1.2}
\lim_{\mathbf{x}\rightarrow\mathbf{x}^{{\star}}}
\frac{\left\|\delta\left(\mathbf{x}\right)\right\|}{
\left\|\mathbf{x}-\mathbf{x}^{{\star}}\right\|}=0
\end{equation}

\item{\textbf{(C2)}}: The step-size sequence, $\left\{\alpha_t\right\}$ is of the form,
\begin{equation}
\label{assc2.1} \alpha_t.=\frac{a}{t+1},~~\mathrm{for\,any}\, t\geq 0,
\end{equation}
where $a>0$ is a constant.

\item{\textbf{(C3)}}: Let $I$ be the $d\times d$ identity matrix and
$a,B$ as in~\eqref{assc2.1} and~\eqref{assc1.1}, respectively. Then, the matrix $\Sigma= aB+\frac{1}{2}I$ is stable.
\item{\textbf{(C4)}}: The entries of the matrices, for any $
t\geq 0, \mathbf{x} \in{R}^{d}$,
\[
\mathbf{A}\left(t,\mathbf{x}\right)=\mathbb{E}\left[{\boldsymbol{\gamma}}\left(t+1,\mathbf{x},\omega\right)
{\boldsymbol{\gamma}}^{\top}\left(t+1,\mathbf{x},\omega\right)\right],
\]
are finite, and the following limit exists: $\lim_{t\rightarrow\infty,~\mathbf{x}\rightarrow\mathbf{x}^{{\star}}} \mathbf{A}\left(t,\mathbf{x}\right) =\mathcal{S}_{0}$
\item{\textbf{(C5)}}: There exists $\epsilon>0$, such that
\begin{equation}
\label{assc5.1}
\hspace{-.1cm}
\lim_{R\rightarrow\infty}\sup_{\left\|\mathbf{x}-\mathbf{x}^{{\star}}\right\|
<\epsilon}\sup_{t\geq 0} \int_{\left\|{\boldsymbol{\gamma}}\left(t+1,\mathbf{x},\omega\right)\right\|>R}
\hspace{-1.4cm}
\left\|{\boldsymbol{\gamma}}\left(t+1,\mathbf{x},\omega\right)\right\|^{2}dP=0
\end{equation}
\end{itemize}
Then we have the following:

Let Assumptions~\textbf{(B1)-(B5)} hold for
$\left\{\mathbf{x}^t\right\}$ in \eqref{RM:1}. Then,
starting from an arbitrary initial state, the process
$\left\{\mathbf{x}^t\right\}$ converges a.s. to
$\mathbf{x}^{{\star}}$.

The normalized process, $\left\{\sqrt{t}\left(\mathbf{x}^t- \mathbf{x}^{{\star}}\right)\right\}$, is asymptotically normal if, besides Assumptions~\textbf{(B1)-(B5)}, Assumptions~\textbf{(C1)-(C5)} are also satisfied. In particular, as $t\rightarrow\infty$, we have:
\begin{equation}
\label{RMres2}
\sqrt{t}\left(\mathbf{x}^t-\mathbf{x}^{{\star}}\right) \xrightarrow{d}. \mathbb{N}({0},\mathcal{S}).
\end{equation}
Also, the asymptotic covariance  $\mathcal{S}$ of the multivariate distribution 
$\mathbb{N}({0},\mathcal{S})$  is
\begin{equation}
\label{RMres3} \mathcal{S}=a^{2}\int_{0}^{\infty}e^{v\,\Sigma }\mathcal{S}_{0}e^{v\,\Sigma^{\top}}dv.
\end{equation}
\end{theorem}
\textit{Proof.} For a proof see \cite{nevelson1976stochastic} (c.f. Theorems 4.4.4, 6.6.1).

We also make use of the following Theorem that is a slight modification of Lemmas 4 and 5 in \cite{SoummyaAppendix}.
\begin{theorem}
\label{lemma-estimation}
Let $z^t$ be a nonnegative (deterministic) sequence satisfying:
\[
z^{t+1} \leq (1-r_1^t)\,z_1^t + r_2^t,
\]
where $\{r_1^t\}$ and $\{r_2^t\}$
are deterministic sequences with

\begin{eqnarray*}
\frac{a_1}{(t+1)^{\delta_1}}
\leq r_1^t \leq 1 \,\,\,\mathrm{and}\,\,\,
r_2^t \leq
\frac{a_2}{(t+1)^{\delta_2}},
\end{eqnarray*}
with $a_1 , a_2 >0,$ and $ \delta_2 > \delta_1>0$.
Then, the following holds: (1) If $\delta_1<1$, then
$z^t = O(\frac{1}{t^{\delta_2-\delta_1}})$;
(2)
If $\delta_1=1$, then
$z^t = O(\frac{1}{t^{\delta_2-1}})$
provided that $a_1>\delta_2-\delta_1$;
(3) if $\delta_1=1$ and
$a_1 < \delta_2-1$, then
 $z^t=O(\frac{1}{t^\zeta})$, 
for any
$\zeta < a_1$.
\end{theorem}

\subsection*{\textbf{B}. A demonstration that the linear SGD's iterate sequence has infinite variance}
\label{section-appendixB}
We provide here a simple demonstration that the linear SGD's iterate sequence 
has infinite variance under the setting of Assumptions~\ref{assumption-f-i}, 
\ref{assumption-nu-comp-wise-iid-noise}, and Assumption~\ref{assumption-psi-comp-wise-iid-noise}, condition~3., holds.

More precisely, assume that 
the gradient noise $\nu^t$ has infinite variance. Consider algorithm~\eqref{eqn-SGD-nonlinear} for solving problem~(1) with $f:\,\mathbb R \mapsto \mathbb R$, $f(x)=\frac{x^2}{2}$, with $\Psi$ being the identity function. Further, consider arbitrary sequence of positive step-sizes $\{\alpha_t\}$. Then, we have:

\begin{equation}
\label{eqn-append-counterex}
x^{t+1} = (1-\alpha_t)\,x^t - \alpha_t\,\nu^t,\,\,t=0,1,...,
\end{equation}
with arbitrary deterministic initialization~$x^0 \in {\mathbb R}$. Then, squaring \eqref{eqn-append-counterex}, using the independence of $x^t$ and $\nu^t$, and the fact that $\nu^t$ has zero mean, we get:
 $
\mathbb{E}\left[(x^{t+1})^2\right] $
$
=
(1-\alpha_t)^2\,\mathbb{E}\left[(x^{t})^2\right] 
$
$+ \alpha_t^2 
\,\mathbb{E}[(\nu^t)^2] $
$
\geq \alpha_t^2\, 
$
$\mathbb{E}[(\nu^t)^2],\,\,t=0,1,...$ 
Taking expectation and using the fact that $\mathbb{E}[(\nu^t)^2]=+\infty$, we see that $\mathbb{E}\left[(x^{t})^2\right] = +\infty$, for any $t \geq 1$.

\subsection*{\textbf{C}. Extension of Theorem~\ref{theorem-MSE-comp-Wise} for gradient noise vector with mutually dependent entries} \label{section-appendixC} 

We show that Theorem~\ref{theorem-MSE-comp-Wise} continues to hold if Assumptions \ref{assumption-nu-comp-wise-iid-noise}, parts 2. and 3., are relaxed, i.e., when we have an i.i.d. zero mean noise vector sequence $\{\nu^t\}$ with a joint pdf $p:\,{\mathbb R}^d \mapsto \mathbb R$. In more detail, we provide an extension of Lemma~\ref{lemma2.2} but for component-wise nonlinearities. Namely, as in Lemma~\ref{lemma2.2}, consider, for a fixed $y \neq 0$:
   \begin{equation}
   \int \psi(\mathbf{y}+\mathbf{u})^\top \mathbf{y} \,p(\mathbf{u})\,d\mathbf{u}.
   \end{equation}
   As, for $\mathbf{a} \in {\mathbb R}^d$, we have $\Psi(\mathbf{a}) = (\Psi(a_1),...,\Psi(a_d))^\top$ (component-wise nonlinearity), we have:
    \begin{eqnarray*}
    &\,& \int \psi(\mathbf{y}+\mathbf{u})^\top 
    \mathbf{y} \,p(\mathbf{u})\,d\mathbf{u} = 
    \int \left( \sum_{i=1}^d \psi(y_i+u_i) y_i \right)\,p(\mathbf{u})\,d\mathbf{u}\\
    &\,& = 
    \sum_{i=1}^d 
    \int \left(  \psi(y_i+u_i) y_i \right)\,p(\mathbf{u})\,d\mathbf{u} = 
    \sum_{i=1}^d \int \left(  \psi(y_i+u_i) y_i \right)\,p_i(u_i)\,du_i,
    \end{eqnarray*}
    where $p_i(u_i)$ is the marginal pdf of the $i$-th component of $\boldsymbol{\nu}^t$. It is easy to show, as $p(\mathbf{u}) = p(-\mathbf{u})$, $\mathbf{u} \in {\mathbb R}^d$, 
    that, for any $i=1,...,d$, we have  $p_i(u) = p_i(-u)$, $u \in {\mathbb R}$. Define $\phi_i(a) - \int \psi(a+u)p_i(u) du$. Note that $\phi_i(a)$ now obeys Lemma~{5.1}. In particular, $\phi_i$ is also  odd, and hence:
    $
    \int \psi(\mathbf{y}+\mathbf{u})^\top $
    $
    \mathbf{u} \,p(\mathbf{u})\,d\mathbf{u} $
    $\geq
    \sum_{i=1}^d  |\phi_i(y_i)|\,|y_i|.
    $ 
    The proof now proceeds analogously to that of Theorem~\ref{theorem-MSE-comp-Wise}.
   
\newpage

\bibliographystyle{unsrt}
\bibliography{references}

\end{document}